# Tractable infinite-dimensional model for long-term environmental impact assessment of long-memory decay processes

Hidekazu Yoshioka[1,] *, Kunihiko Hamagami[2]

[1] Graduate School of Advanced Science and Technology, Japan Advanced Institute of Science and Technology, 1-1 Asahidai, Nomi, Ishikawa, Japan

[2] Faculty of Agriculture, Department of Food Production and Environmental Management, Iwate University, 3-18-8 Ueda, Morioka, Japan

* Corresponding author: yoshih@jaist.ac.jp, ORCID: 0000-0002-5293-3246

***Abstract***

Focusing on the assessment of benthic algae blooms that decay subexponentially, we propose a tractable (solvable in a semi-explicit form) and well-defined (that does not diverge) environmental index for the impact assessment of long-memory processes under model uncertainties. Our target system generates long-memory decay through an infinite superposition of multiscale processes. The sensitivity of the environmental index can be controlled by the degree of model uncertainty in terms of the relative entropy and nonexponential discount; hence, we apply a long-memory discount to evaluate long-memory processes. In our framework, the evaluation of the environmental index is reduced to finding a proper solution to an infinite-dimensional extended Hamilton–Jacobi–Bellman system. We solve this system, verify a reasonable solution, and numerically handle it by using a quantization technique. Finally, we present a demonstrative application of the proposed framework to benthic algae population dynamics in river environments based on a laboratory experiment. This paper offers a tractable framework for the assessment of persistent environmental phenomena.

***Keywords***

Long-memory decay; Environmental index; Time-inconsistent optimal control; Nonexponential discount; Benthic algae population dynamics

***Statements & declarations***

**Acknowledgments:** We are deeply grateful for the valuable comments from the anonymous reviewers.

**Funding:** This study was supported by the Japan Society for the Promotion of Science (KAKENHI No. 25K07931) and the Japan Science and Technology Agency (PRESTO No. JPMJPR24KE).

**Conflict of Interests:** The authors have no relevant financial or nonfinancial interests to disclose.

**Data Availability:** The data will be made available upon reasonable request to the corresponding author.

**Declaration of Generative AI in Scientific Writing:** The authors did not use generative AI to obtain the results of this study.

# 1. Introduction

## 1.1 Background and motivation

Long-memory processes are time-dependent variables with persistent memory that decays subexponentially. These processes are ubiquitous in our world; examples include but are not limited to precipitation (Fritsch et al., 2025)[1], socioenvironmental droughts (Granata and Di Nunno, 2026)[2], air pollution (Shi et al., 2025)[3], algae blooms (Lui et al., 2007)[4], transport phenomena (Qu et al., 2026)[5], clean energy power trends (Liu et al., 2025)[6], electricity prices (Balagula et al., 2025)[7], and natural language (Stanisz et al., 2024)[8]. Therefore, impact assessment of long-term memory processes is pivotal for predicting their influence on human societies as well as local and global environments.

Mathematically, a wide class of long-memory processes has representation formulae in which a process is described as a superposition (i.e., summation or integration) of short-memory (exponential-memory) processes. A classical example is fractional Brownian motion, which is expressed as a stochastic integration of infinitely many correlated Ornstein–Uhlenbeck processes (Muravlev, 2011)[9]. Other examples include superpositions of independent Ornstein–Uhlenbeck processes (Barndorff-Nielsen, 2001; Fasen and Klüppelberg, 2005)[10,11]. These mathematical techniques, known as Markovian embedding and Markovian lifts (we use the latter in this paper), have also been applied to more complex models, such as the generalized Langevin system (Siegle et al., 2011)[12], rough volatility models (Abi Jaber and El Euch, 2019)[13], sandwiched Volterra processes (Di Nunno and Yurchenko-Tytarenko, 2026)[14], gray Brownian motions (Bock et al., 2020)[15], fractional Brownian bridges (Nobis et al., 2025)[16], American option models (Aichinger and Desmettre, 2025)[17], nonlinear Hawkes processes (Kanazawa and Dechant, 2025)[18], and jump-driven water quantity-quality models (Yoshioka and Yoshioka, 2025a; Yoshioka and Yoshioka, 2025b)[19,20]. A pivotal advantage of Markovian lifts is that the processes to be superposed are often formally simpler than the original long-memory processes, which have been effectively used in the numerical approximation of long-memory processes based on quantization techniques for simulation and assessment (Alfonsi, 2025; Dall'Acqua et al., 2024; Hainaut et al., 2025; Lin et al., 2025; Yoshioka et al., 2024)[21,22,23,24,25]. Markovian lifts are also connected to superstatistics, which statistically generate a complex system by superposing simpler lifts (Davis et al., 2025; Sánchez, 2025)[26,27].

Model uncertainties, namely, misspecifications of parameters and coefficients in a mathematical model due to data limitations, are crucial in operating a stochastic process model similar to a long-memory process; examples include regional precipitation (Moges et al., 2022)[28], air quality and weather indices (Huang et al., 2025)[29], financial indices with heavy tails (Bergmann and Oliveira et al., 2026)[30], market volatility (Lehrer et al., 2021)[31], river discharge (Yoshioka and Yoshioka, 2024a)[32] and water quality (Yoshioka and Yoshioka, 2024b)[33]. A well-established and analytically tractable framework for modeling stochastic processes subject to uncertainties is the optimal control methodology, in which uncertainties are expressed by control variables to maximize a disutility index (Jaimungal et al., 2024)[34]. Insurance models within this framework have recently been studied (Kroell et al., 2024; Kroell et al., 2025)[35,36]. The cost of control is given by the relative entropy between the benchmark (or estimated) and true models, which has useful convexity that poses a control problem and allows for computation in complex applications (Horiguchi and Kobayashi, 2025)[37]. The control approach for model uncertainties has not been extensively studied considering the long-term environmental impact assessments discussed below.

A common characteristic of long-memory processes is their persistence, such that their influence decays slowly, which can be briefly explained through the following two relationships:

$$\int_0^{+\infty} e^{-t} \mathrm{d}t = 1 < +\infty , \tag{1}$$

while

$$\int_0^{+\infty} \left( \frac{1}{1+t} \right)^{\alpha} \mathrm{d}t = 1 < +\infty \quad \text{if and only if} \quad \alpha > 1 . \tag{2}$$

The first relationship (1) simply states that the cumulative influence of exponential decay is bounded, whereas the second relationship (2) shows that the cumulative influence of subexponential decay is unbounded if $0 < \alpha \leq 1$, that is, the truly long-memory case. Yoshioka and Hamagami (2024)[38] discussed that the population of nuisance riverine benthic algae decays subexponentially during flood events, suggesting their persistent influence on aquatic environments and ecosystems. The population decay of benthic algae is a key process in their population dynamics, which dominate slow water flows (Haddadchi et al., 2020)[39]. Moreover, long-term prediction of benthic algae blooms is difficult, suggesting the importance of finding an effective method for environmental assessment (Luce et al., 2010)[40].

A possible solution to derive a bounded integral of a long-memory process is to incorporate an exponential discount rate, as commonly used in stochastic calculus and control (e.g., Remark 3.2.2 in Pham (2009)[41]), which in our context is given as follows:

$$\int_0^{+\infty} e^{-\delta t} \left( \frac{1}{1+t} \right)^{\alpha} \mathrm{d}t \leq \int_0^{+\infty} e^{-\delta t} \mathrm{d}t = \frac{1}{\delta} < +\infty \quad \text{for any} \quad \delta > 0 . \tag{3}$$

The discount rate $\delta > 0$ measures how future states will be weighted in the integral, and a larger value of $\delta$ implies a more myopic integrand. However, this methodology is excessively strong considering the sensitivity of the integral to the key parameter $\alpha$, which represents subexponential decay, because $e^{-\delta t}$ decreases significantly faster than $\left( \frac{1}{1+t} \right)^{\alpha}$. This implies that the integral in (3) is mainly controlled by the decay but not the long-memory process, with the latter being the target. This finding motivates us to consider a nonexponential decay, such that the target integral is still bounded, but the decay is not excessively fast compared with the long-memory process. This theoretical issue would be more pronounced in applications where one needs to evaluate a worst-case integral index that is a version of (3) maximized with respect to some uncertain model parameter and/or coefficients (e.g., Jaimungal et al., 2024)[34]. Therefore, the traditional approach based on the exponential discount (Babonneau et al., 2025; Kafash and Nikooeinejad, 2025; Lykina et al., 2022; Marsiglio and Masoudi, 2025; Yoshioka et al., 2020)[42-46] can be improved via a nonexponential alternative unless the target process exhibits exponential decay, which is the motivation of this study.

Consequently, a candidate for a long-term environmental index is given by the integration of a long-memory process that is maximized subject to uncertainties and nonexponential discounts. This is a difficult topic compared with classical control models because the dynamic programming principle, which is a powerful tool for stochastic control problems, cannot be employed owing to the nonexponential discount; however, in economics, significant progress has been made in resolving control problems with nonexponential discounts, as reviewed by Björk et al. (2021)[47]: the time-inconsistent control theory.

A formal difference between classical and time-inconsistent control theories is that optimal control is attained by solving the Hamilton–Jacobi–Bellman (HJB) equation in the former and an extended (i.e., more complicated) HJB equation system in the latter. In the context of time-inconsistent control theory, control problems subject to nonexponential discounts have been effectively solved by focusing on tractable utility maximization (Ekeland and Lazrak, 2010)[48], consumption and investment (Chen et al., 2024; Liu et al., 2020)[49,50], contracts (Cetemen et al., 2023)[51], dividends (Li et al., 2016; Zhu et al., 2020)[52,53], portfolios (Kang et al., 2025; Mbodji and Pirvu, 2025)[54,55], and optimal stopping problems (Ebert et al., 2020)[56]. However, this time-inconsistent formulation has not yet been applied to long-term environmental impact assessments.

### 1.2 Aim and contribution

The aim of this study is to formulate an effective index for long-term environmental impact assessment of long-memory processes. In this paper, by a long-memory process, we mean a time-dependent variable that decays slower than exponentially, and we do not focus on subexponential growth. The most compelling motivation in this paper is to propose a theory where persistent decay, slow discount, and model ambiguity interact.

The main goal of this study is to generalize the following relationship between a more general discount factor and long-memory processes:

$$\int_0^{+\infty} \underbrace{\left(\frac{1}{1+at}\right)^{\beta}}_{\text{Discount}} \underbrace{\left(\frac{1}{1+t}\right)^{\alpha}}_{\text{Long-memory process}} \mathrm{d}t < +\infty \quad \text{for any} \quad a > 0 \quad \text{if and only if} \quad \alpha + \beta > 1 . \tag{4}$$

Specifically, the left-hand side of (4) is bounded if and only if $\beta > \max\{0, 1-\alpha\}$, which holds true even if $\alpha \in (0,1]$, with which (2) fails. In particular, the nonexponential discount is weaker than the exponential discount (i.e., (3)) because of the significantly slower decay.

The process we consider is of the superposition type studied by Yoshioka and Hamagami (2025)[57] as a representative ecological model of subexponential decay, which was motivated by the subexponential decay of the benthic algae population on the riverbed (**Figure 1**). This model is a superposition of infinitely many two-state Markov chains and has an approximate finite-dimensional representation. Superposition emerges as a deterministic system, which in our case is a continuum of ordinary differential equations. This is consistent with time-inconsistent control models with large problem dimensions (Liang and Zhang, 2024; Niu and Zou, 2024)[58,59] and the control-theoretic approach (e.g.,

Jaimungal et al., 2024)[34] for model uncertainties. Similar methodologies have been applied to time-consistent versions (Yoshioka and Yoshioka, 2023)[60]. The algae cover, as shown on the right side of **Figure 1**, has been found to decay subexponentially over time, and its population dynamics modeling needs to be considered (Yoshioka and Hamagami, 2025)[57].

We assume that model uncertainties appear in the memory length of the target process, which is a key quantity; estimating the memory length without errors in applications would be difficult. Under model uncertainties, we present an environmental index subject to a nonexponential discount. We demonstrate that it can be evaluated as a semi-explicit solution to an infinite-dimensional extended HJB system whose optimality is proven under certain conditions. Infinite-dimensional systems and their controls have been discussed theoretically (e.g., Boucekkine et al., 2013; Kanazawa and Sornette, 2024; Li and Yong, 1995)[61,62,63], but their implementation in environmental problems is still scarce. The literature suggests that an infinite-dimensional system arises as a large-size limit of a finite-dimensional system, which in our case is a quantization version (Flandoli et al., 2025; Yoshioka and Hamagami, 2025; Yu et al., 2021)[57,64,65] that is computationally implementable. From this perspective, we can effectively approximate an infinite-dimensional system via finite-dimensional systems via a quantization technique. Note that Galerkin approximation methods would be an option for approximating infinite-dimensional stochastic models (Karbach, 2026)[66], while they are not able to preserve signs of the system solutions; our approximation is simpler, more intrinsic, and does not encounter this issue.

We present a demonstrative application of the proposed framework for long-term environmental impact assessment of the population dynamics of nuisance benthic algae in **Figure 1** (see **Section 4** for further explanations). Model uncertainties in this case arise from the limitation of data quantity, with which the timescale of population decay is difficult to estimate without errors. Long-term algae blooms in water bodies are a critical environmental issue (Igwaran et al., 2024; Sun et al., 2023; Xu et al., 2025; Yang et al., 2025)[67,68,69,70]. Hence, exploring an effective environmental index for their assessment is considered important. We computationally investigate the environmental impacts of the persistent algae population and the sensitivity of the corresponding environmental index with a nonexponential discount and model uncertainties. Consequently, this study contributes to the formulation, analysis, and application of a new mathematical framework for impact assessment of long-memory processes.

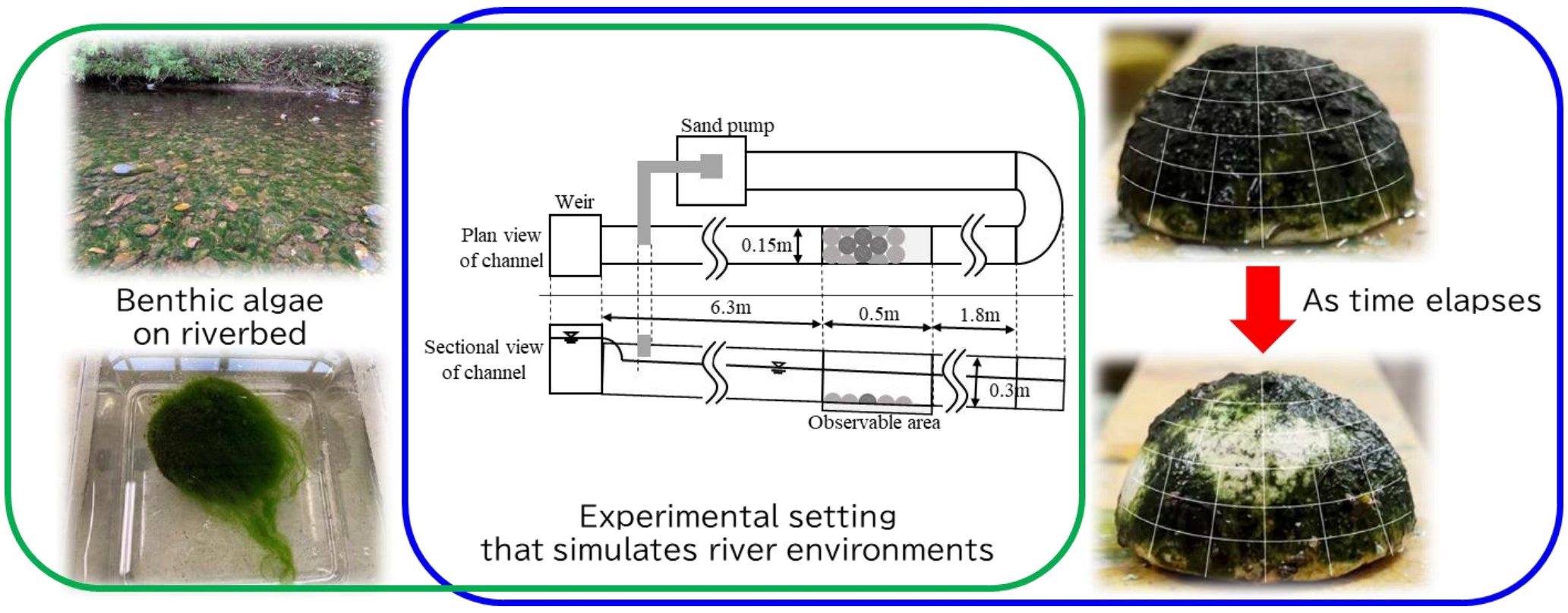

**Figure 1.** Conceptual image of the motivational problem of this study.

### 1.3 Structure of this paper

The rest of this paper is organized as follows. **Section 2** introduces the proposed long-memory process model and its finite-dimensional version. **Section 3** presents the environmental index and discusses the associated extended HJB system. **Section 4** demonstrates an application of the proposed mathematical framework. **Section 5** concludes the paper. **Appendix** presents the proofs and auxiliary results.

## 2. Long-memory process

### 2.1 Finite-dimensional version

First, a finite-dimensional version of the proposed model is discussed to explain how the system and its model uncertainties arise.

#### 2.1.1 Model without uncertainties

Our dynamics originate from a system of stochastic differential equations (SDEs). The explanation below focuses on an application to the benthic algae population but can be extended to other population dynamics with suitable adaptations.

We follow the problem setting in Yoshioka and Hamagami (2025)[57]. Consider a habitat $\Omega$ with a finite area $|\Omega| > 0$, such as a river reach (e.g., left side of **Figure 1**). Without any loss of generality, we set $|\Omega| = 1$. We consider the continuous-time population dynamics of a single species living in $\Omega$. The population changes with time according to decay and growth, and we assume that the growth rate is constant over $\Omega$ but that the decay rate varies. In the context of benthic algae population dynamics, this setting corresponds to the situation in which growth by photosynthesis and decay is caused by abrasion, which has been found to be heterogeneous over a riverbed that is typically nonflat (Graba et al., 2010; Hamagami et al., 2024)[71,72]. We divide $\Omega$ into $N \in \mathbb{N}$ subdomains $\{\Omega_i\}_{i=1,2,3,\ldots,N}$ on which the decay rate of the population is approximated to be constant. The area of $\Omega_i$ is $|\Omega_i|$, and we have $|\Omega| = \sum_{i=1}^{N} |\Omega_i|$. For modeling simplicity, the population in each $\Omega_i$ is either 1 (a population exists) or 0 (otherwise). Time is denoted as $t \geq 0$, and the population density in $\Omega_i$ at time $t$ is denoted as $x_{i,t}$.

Our model starts with the following system of SDEs:

$$\underbrace{\mathrm{d}x_{i,t}}_{\text{Increment}} = \underbrace{-x_{i,t-}\mathrm{d}J_{i,t}}_{\text{Decay}} + \underbrace{\left(1 - x_{i,t-}\right)\mathrm{d}K_{i,t}}_{\text{Growth}}, \quad t > 0 \tag{5}$$

with initial conditions $\{x_{i,0}\}_{i=1,2,3,\ldots,N} \in \{0,1\}^N$, where $\{J_{\cdot,i}\}_{i=1,2,3,\ldots,N}$ is a collection of mutually independent Poisson processes; the jump rate (i.e., decay rate) of $J_{\cdot,i}$ is a constant $r_i > 0$, which is the decay rate of algae in $\Omega_i$; and $\{K_{\cdot,i}\}_{i=1,2,3,\ldots,N}$ is another collection of mutually independent Poisson processes with a common jump rate $R > 0$. We also assume that the elements of $\{J_{\bullet,i}\}_{i=1,2,3,\ldots,N}$ and $\{K_{\bullet,i}\}_{i=1,2,3,\ldots,N}$ are mutually independent. The population $X_t^{(N)}$ in $\Omega$ at time $t$ is the sum

$$X_t^{(N)} = \sum_{i=1}^{N} |\Omega_i| x_{i,t}, \quad t \geq 0. \tag{6}$$

Note here that the growth rate of each $x_{i,t}$ is exponential with the common rate $R > 0$.

### 2.1.2 Model uncertainties

The model uncertainties in this study refer to misspecifications in the decay rates $\{r_i\}_{i=1,2,3,\ldots,N}$, which control the decay, namely, memory in our context, of the population (see **Remark 1** presented later). Our formulation is based on the model uncertainties for jump processes (Jaimungal et al., 2024)[34]. We assume that uncertainties exist such that $\{r_i\}_{i=1,2,3,\ldots,N}$ (estimated model) is distorted to $\{\phi_i r_i\}_{i=1,2,3,\ldots,N}$ (i.e., possibly true model) with positive sequences $\{\phi_i\}_{i=1,2,3,\ldots,N}$ such that for $\Omega_i$, misspecifications are present if $\phi_i \neq 1$, whereas none exist if $\phi_i \equiv 1$.

Under the misspecification setting, the SDE (5) becomes

$$\mathrm{d}x_{i,t} = -x_{i,t-}\mathrm{d}\hat{J}_{i,t} + \left(1 - x_{i,t-}\right)\mathrm{d}K_{i,t}, \quad t > 0, \tag{7}$$

where $\{\hat{J}_{\bullet,i}\}_{i=1,2,3,\ldots,N}$ is a collection of mutually independent Poisson processes, and the jump rate of $\hat{J}_{\bullet,i}$ is $\phi_i r_i > 0$. The unit-time relative entropy between the benchmark and distorted models for each $i$ is $\phi_i \ln\phi_i - \phi_i + 1$, which is a nonnegative convex function of $\phi_i > 0$ that is minimized at $\phi_i = 1$, with a minimum value of 0. This unit-time relative is a building block of the environmental index.

## 2.2 Finite- to infinite-dimensions

We consider an infinite-dimensional version of system (7), which can be considered a generalization of the finite-dimensional system. We can pass the limit $N \to +\infty$ in (7) as performed in Proposition 1 in Yoshioka and Hamagami (2025)[57] because the stochastic system (7) is linear with respect to each $x_i$. We assume that each $r_i$ is generated by a probability measure $p$ of positive random variables, which decays exponentially fast for $r \to +\infty$, such that $|\Omega_i|$ can be approximately understood as the probability of selecting $r_i$. In the framework of Yoshioka and Hamagami (2025)[57], the population at time $t$ (a.s. convergence limit of $X_t^{(N)}$) is given by

$$X_t = \int_0^{+\infty} x_t(r)\, p(\mathrm{d}r), \quad t \geq 0 \tag{8}$$

with

$$\underbrace{\frac{\partial x_t(r)}{\partial t}}_{\text{Variation}} = \underbrace{-r\phi_t(r)x_t(r)}_{\text{Decay}} + \underbrace{R(1-x_t(r))}_{\text{Growth}}, \quad t>0 \text{ and } r>0 \text{ for } p \text{ a.e.} \tag{9}$$

subject to an initial condition $x_0$, which is a continuous function in $(0,+\infty)$ that is uniformly bounded between 0 and 1. Here, $\phi$ is a positive function in $(0,+\infty)$ that represents model uncertainties. Model uncertainties arise in the first term on the right-hand side of (9). Here, $x_t(r)$ would be interpreted as an expected occupation fraction for $r$ at $t$.

The equation (9) can be solved as follows:

$$x_t(r) = x_0(r)e^{-\int_0^t (r\phi_s(r)+R)\mathrm{d}s} + R\int_0^t e^{-\int_s^t (r\phi_u(r)+R)\mathrm{d}u}\mathrm{d}s \tag{10}$$

and hence

$$X_t = \int_0^{+\infty}\left\{x_0(r)e^{-\int_0^t (r\phi_s(r)+R)\mathrm{d}s} + R\int_0^t e^{-\int_s^t (r\phi_u(r)+R)\mathrm{d}u}\mathrm{d}s\right\}p(\mathrm{d}r). \tag{11}$$

In this setting, the unit-time relative entropy for each $r>0$ is given by the continuum version $\phi_t(r)\ln\phi_t(r)-\phi_t(r)+1$. The solution (10) satisfies (9) pointwise and is bounded between 0 and 1 globally in time.

***Remark 1*** The long-memory nature of the process $X$, particularly its decay, originates from the singularity of $p$ near $r=0$. Indeed, if $p$ has a density that behaves asymptotically $r^{\alpha-1}$ with $\alpha>0$, then

$$\int_0^{+\infty} e^{-rt}p(\mathrm{d}r) = O(t^{-\alpha}) \quad \text{for } t>>1, \tag{12}$$

where we assume that $x_0(r)\equiv 1$, which corresponds to the case without growth ($R=0$). The long-memory nature of our system thus arises from heterogeneous population decay.

We conclude **Section 2** by defining a class of bounded functions that will be used in the sequel. The collection of $p$ a.e. (almost everywhere) measurable functions $f:(0,+\infty)\to\mathbb{R}$ bounded between 0 and 1 is denoted by $B(0,+\infty)$. If $x_0\in B(0,+\infty)$, then $x_t\in B(0,+\infty)$ for any $t>0$ by (10). Notably, the space $E=L^\infty(p)$ of $p$ a.e. measurable and bounded functions equipped with the maximum norm, is a Banach space (Chapter 3.3 in Clason (2020)[73]).

## 3. Environmental index

In the main text, we consider an environmental index for the infinite-dimensional system. See **Appendix** for the environmental index in a finite-dimensional setting.

### 3.1 Control-theoretic formulation

The objective function is as follows:

$$\theta(x_0;\phi)=\int_0^{+\infty}\left(\underbrace{\int_0^{+\infty}e^{-\delta s}X_s\mathrm{d}s}_{\text{Cumulative disutility}}\underbrace{-\int_0^{+\infty}e^{-\delta s}\left(\int_0^{+\infty}\frac{r}{\eta(r)}x_s(r)\left(\phi_s(r)\ln\phi_s(r)-\phi_s(r)+1\right)p(\mathrm{d}r)\right)\mathrm{d}s}_{\text{Uncertainty penalization}}\right)\mu(\mathrm{d}\delta), \quad (13)$$

where $\mu$ is the probability measure of the distributed discount rate $\delta>0$, which does not have a Dirac mass at $\delta=0$. This index $\theta$ depends on the initial condition $x_0$ and uncertainty $\phi$ and contains two terms; the first term on the right-hand side of (13) represents cumulative disutility caused by the existence of the population, and the second term represents the penalization of uncertainty whose size is measured through the convex function $\phi_s(r)\ln\phi_s(r)-\phi_s(r)+1$ of $\phi_s(r)$, which comes from the relative entropy between the benchmark ($\phi\equiv 1$) and distorted models ($\phi\neq 1$) (Jaimungal et al., 2024)[34]. In this case, the factor $x_s(r)$ implies that the uncertainty for $r$ needs to be accounted for only when a corresponding population ($x_s(r)>0$) is present. Another factor, $\frac{r}{\eta(r)}$, which serves as a weighting coefficient between the two terms on the right-hand side of (13), measures the significance of the uncertainties in $\theta$. More specifically, a larger $\eta$ corresponds to a more uncertainty-averse (pessimistic) index $\theta$. The functional form of $\eta$ is critical for the well-posedness of our environmental index (**Proposition 1** presented later).

Finally, the discount factor is incorporated into (13) through distributed discounted rates via the probability measure $\mu$. Its regularity near $\delta=0$ determines the nonexponential nature of the discount, as in the case of $p$ (**Remark 1**). By allowing for the exchange of parameters in integrals, we can rewrite (13) as

$$\theta(x_0;\phi)=\int_0^{+\infty}\Delta_s X_s\mathrm{d}s-\int_0^{+\infty}\Delta_s\int_0^{+\infty}\frac{r}{\eta(r)}x_s(r)\left(\phi_s(r)\ln\phi_s(r)-\phi_s(r)+1\right)p(\mathrm{d}r)\mathrm{d}s \quad (14)$$

with a discount factor that is not necessarily exponential, which is another superposition (of exponential functions) employed in the present mathematical framework:

$$\Delta_t=\int_0^{+\infty}e^{-\delta t}\mu(\mathrm{d}\delta),\ t>0, \quad (15)$$

where $\lim_{t\to+\infty}\Delta_t=0$ because $\mu$ does not have a Dirac mass at $\delta=0$.

The index $\theta$ is maximized with respect to $\phi$. The admissible set $\mathbb{F}$ of controls $\phi$ is a collection of positive and measurable functions in $(0,+\infty)$ such that system (9) admits a unique bounded solution (10) that is continuously differentiable in time, satisfying

$$\int_0^{+\infty}\int_0^{+\infty}e^{-\delta s}X_s\mathrm{d}s\mu(\mathrm{d}\delta)<+\infty \quad (16)$$

and

$$\int_0^{+\infty}\int_0^{+\infty}e^{-\delta s}\left(\int_0^{+\infty}\frac{r}{\eta(r)}x_s(r)\left(\phi_s(r)\ln\phi_s(r)-\phi_s(r)+1\right)p(\mathrm{d}r)\right)\mathrm{d}s\mu(\mathrm{d}\delta)<+\infty . \quad (17)$$

The value function, which is the maximized $\theta$ with respect to $\phi$, is set to

$$V(x_0) = \text{"}\sup_{\phi\in\mathbb{F}}\text{"}\,\theta(x_0;\phi). \tag{18}$$

Here, "sup" is taken with respect to all $\phi \in \mathbb{F}$ that satisfy the conditions of "Markovian" equilibrium controls described in **Definition 1** (e.g., Definition 15.3 in Björk et al. (2021)[47]).

***Definition 1***

*Consider a Markovian control law $\hat{\phi} \in \mathbb{F}$. Set an arbitrary Markovian control $\phi \in \mathbb{F}$. Define $\phi_\varepsilon \in \mathbb{F}$ and a fixed real number $\delta > 0$ such that*

$$\phi_{\varepsilon,s}(r,y) = \begin{cases} \phi_s(r,y) & (0 \le s < \varepsilon,\ y \in E) \\ \hat{\phi}(r,y) & (\varepsilon \le s,\ y \in E) \end{cases},\quad r > 0 \ \text{ for } \ p \ \text{ a.e.} \tag{19}$$

*If*

$$\liminf_{\varepsilon\to+0} \frac{\theta(x;\hat{\phi}) - \theta(x;\phi_\varepsilon)}{\varepsilon} \ge 0 \tag{20}$$

*for all $\phi \in \mathbb{F}$ and $x \in B(0,+\infty)$, then we call $\hat{\phi}$ an equilibrium control and $\theta(\cdot;\hat{\phi})$ the value function $V$.*

When there is no confusion, we omit the second argument in $\phi$. The theoretical goal of this study is to determine the value function $V$, which is the environmental index used herein. This task is achieved by solving an extended HJB system, as explained in the following section.

***Remark 2*** The admissibility condition of $\phi$, particularly in (16) and (17), guarantees that each term in (13) is bounded. Without these, a possibility exists that (16) and (17) are unbounded, but their difference is bounded, which is an unphysical situation.

### 3.2 Extended HJB system

In the rest of this paper, we assume

$$\rho_t \equiv \frac{\mathrm{d}}{\mathrm{d}t}\Delta_t = -\int_0^{+\infty} \delta e^{-\delta s}\mu(\mathrm{d}\delta) > -\infty, \tag{21}$$

implying that the probability measure $\mu$ has an average, which is considered not restrictive in applications because we focus on cases where the discount is small.

Our extended HJB system is as follows: for any $f \in E$,

$$\sup_{\phi>0}\left\{ \begin{aligned} &\int_0^{+\infty}\left\{-r\phi(r)f(r) + R(1-f(r))\right\}\Phi_f(r,f)\,p(\mathrm{d}r) \\ &+\int_0^{+\infty} f(r)\,p(\mathrm{d}r) - \int_0^{+\infty}\frac{r}{\eta(r)} f(r)\big(\phi(r)\ln\phi(r) - \phi(r) + 1\big)p(\mathrm{d}r) + \Psi_0(f) \end{aligned}\right\} = 0 \tag{22}$$

coupled with

$$\frac{\partial \Psi_t(f)}{\partial t}+\int_0^{+\infty}\left\{-r\phi^*(r)f(r)+R\left(1-f(r)\right)\right\}\Psi_{t,f}(r,f)p(\mathrm{d}r)$$
$$+\rho_t\left\{\int_0^{+\infty}f(r)p(\mathrm{d}r)-\int_0^{+\infty}\frac{r}{\eta(r)}f(r)\left(\phi^*(r)\ln\phi^*(r)-\phi^*(r)+1\right)p(\mathrm{d}r)\right\}=0 \tag{23}$$

along with the terminal condition (a transversivity criterion to select a reasonable solution that decays at infinity, e.g., )

$$\lim_{t\to+\infty}\Psi_t(f)=0\,. \tag{24}$$

Here, we set

$$\phi^*(\cdot)=\underset{\phi>0}{\arg\max}\left\{\begin{array}{l}-\int_0^{+\infty}r\phi(r)f(r)\Phi_f(r,f)p(\mathrm{d}r)\\ -\int_0^{+\infty}\frac{r}{\eta(r)}f(r)\left(\phi(r)\ln\phi(r)-\phi(r)+1\right)p(\mathrm{d}r)\end{array}\right\}=\exp\left(-\eta(\cdot)\Phi_f(\cdot,f)\right), \tag{25}$$

and the variational (Gateaux) derivative $\Phi_f$ is determined through (p.44 in Li and Yong (1994)[63])

$$\lim_{\varepsilon\to 0}\frac{\Phi(r,f+\varepsilon g)-\Phi(r,f)}{\varepsilon}=\int_0^{+\infty}g(r)\Phi_f(r,f)p(\mathrm{d}r)\quad\text{for any }\ f\in E\,, \tag{26}$$

and that for $\Psi$ is defined similarly.

The reason for the form of the system (22)-(23) becomes clearer in **Proof of Proposition 1 in Appendix**, where this system turns out to arise from the optimality principle (20).

### 3.3 Proper solution

We used a guessed-solution technique with which the extended HJB system can be solved in a semi-explicit form. We assume that the following solution is differentiable in the sense that variational derivatives (e.g., (26)) exist:

$$\Phi(f)=\int_0^{+\infty}A(r)f(r)p(\mathrm{d}r)+B \tag{27}$$

and

$$\Psi_t(f)=\int_0^{+\infty}\left\{\int_0^{+\infty}C(r,\delta)f(r)p(\mathrm{d}r)+D(\delta)\right\}e^{-\delta t}\mu(\mathrm{d}\delta)\,. \tag{28}$$

Here, $A$ is a function in $(0,+\infty)$, $B$ is a constant, $C$ is a function in $(0,+\infty)^2$, and $D$ is a function in $(0,+\infty)$. We assume the boundedness condition

$$\int_0^{+\infty}|A(r)|f(r)p(\mathrm{d}r),\ |B|,\ \int_0^{+\infty}\int_0^{+\infty}|C(r,\delta)|e^{-\delta t}p(\mathrm{d}r)\mu(\mathrm{d}\delta),\ \int_0^{+\infty}|D(\delta)|\mu(\mathrm{d}\delta)<+\infty\,, \tag{29}$$

with which each term in (27) and (28) exists.

The following proposition is our main theoretical result, stating that the extended HJB system admits a semi-explicit solution that corresponds to the environmental index, the value function.

***Proposition 1***

*For* $p$ *a.e., assume*

$$0 < \frac{\eta(r)}{r} < \bar{c} < 1 \quad \textit{for any} \quad r > 0 \tag{30}$$

*with a constant* $\bar{c} \in (0,1)$,

$$\int_0^{+\infty}\int_0^{+\infty} \frac{1}{\delta + r}\mu(\mathrm{d}\delta)p(\mathrm{d}r) < +\infty \quad \textit{if} \quad R = 0, \tag{31}$$

$$\int_0^{+\infty} \frac{1}{\delta}\mu(\mathrm{d}\delta) < +\infty \quad \textit{if} \quad R > 0. \tag{32}$$

*Then, the extended HJB system (22)-(25) admits a semi-explicit solution of the form (27) and (28) where the coefficients* $B, C, D$ *are given as follows: for* $p \times \mu$ *a.e.,*

$$B = R\int_0^{+\infty}\left(1 - \frac{r}{\eta(r)}\left(1 - \left(1 + \eta(r)A(r)\right)e^{-\eta(r)A(r)}\right)\right)\left(\int_0^{+\infty} \frac{1}{\delta\left(\delta + R + re^{-\eta(r)A(r)}\right)}\mu(\mathrm{d}\delta)\right)p(\mathrm{d}r), \tag{33}$$

$$C(r,\delta) = -\frac{\delta}{\delta + re^{-\eta(r)A(r)} + R}\left(1 - \frac{r}{\eta(r)}\left(1 - \left(1 + \eta(r)A(r)\right)e^{-\eta(r)A(r)}\right)\right), \tag{34}$$

*and*

$$D(\delta) = \frac{R}{\delta}\int_0^{+\infty} C(r,\delta)p(\mathrm{d}r). \tag{35}$$

*The coefficient* $A$ *is a unique positive solution* $p$ *a.e. obtained from the following functional equation: for* $p$ *a.e.,*

$$A(r) = \int_0^{+\infty} \frac{1}{\delta + re^{-\eta(r)A(r)} + R}\mu(\mathrm{d}\delta)\left(1 - r\frac{1 - \left(1 + \eta(r)A(r)\right)e^{-\eta(r)A(r)}}{\eta(r)}\right). \tag{36}$$

*Moreover, for* $p$ *a.e., this* $A$ *satisfies the following bound:*

$$0 < A(r) < \bar{A}(r) \equiv -\frac{1}{\eta(r)}\ln\left(1 - \frac{\eta(r)}{r}\right). \tag{37}$$

*Finally, the equilibrium control is given by*

$$\hat{\phi}(r, x_t) = \phi^*(r, x_t) = e^{-\eta(r)A(r)} \tag{38}$$

*and it follows that*

$$V(x_0) = \Phi(x_0) \quad \textit{for any} \quad x_0 \in \mathbb{F}. \tag{39}$$

The conditions required for **Proposition 1** are seemingly complex, but our application study serves as a simpler example where $\eta, p, \mu$ are specified.

***Remark 3*** We briefly discuss the assumptions made in **Proposition 1**. The inequality (30) means that the uncertainty aversion needs to be sufficiently large to define the environmental index well. For that purpose, we consider when the conditions (31)–(32) are satisfied for a tractable case where $p, \mu$ are given by the following gamma distributions, which are representative models for long-memory in the context of superposition processes (e.g., Fasen and Klüppelberg, 2005; Yoshioka and Yoshioka, 2024a)[11,32]:

$$p(\mathrm{d}r) = \frac{1}{\beta_r^{\alpha_r - 1}\Gamma(\alpha_r)} r^{\alpha_r - 1} \mathrm{e}^{-\frac{r}{\beta_r}} \mathrm{d}r \quad \text{and} \quad \mu(\mathrm{d}\delta) = \frac{1}{\beta_\delta^{\alpha_\delta - 1}\Gamma(\alpha_\delta)} \delta^{\alpha_\delta - 1} \mathrm{e}^{-\frac{\delta}{\beta_\delta}} \mathrm{d}\delta, \tag{40}$$

where $\Gamma$ is the Gamma function and $\alpha_r, \beta_r, \alpha_\delta, \beta_\delta > 0$. Then, condition (32) is rewritten as $\alpha_\delta > 1$, and condition (31) is rewritten as $\alpha_r + \alpha_\delta > 1$; condition (31) is satisfied if

$$\int_{\delta^2 + r^2 \leq 1} \frac{r^{\alpha_r - 1}\delta^{\alpha_\delta - 1}}{\delta + r} \mathrm{d}\delta \mathrm{d}r < +\infty \tag{41}$$

because

$$\int_{\delta^2 + r^2 \leq 1} \frac{r^{\alpha_r - 1}\delta^{\alpha_\delta - 1}}{\delta + r} \mathrm{d}\delta \mathrm{d}r = \int_0^1 z^{\alpha_r + \alpha_\delta - 2} \mathrm{d}z \times \underbrace{\int_0^{\frac{\pi}{2}} \frac{(\sin(y))^{\alpha_r - 1}(\cos(y))^{\alpha_\delta - 1}}{\sin(y) + \cos(y)} \mathrm{d}y}_{\text{Bounded}} < +\infty, \tag{42}$$

which holds true when $\alpha_r + \alpha_\delta > 1$. This condition corresponds to that for (4) in **Section 1.**

***Remark 4*** The worst-case model uncertainty (38) is independent of the population but depends on the decay, which is implementable even if we cannot continuously observe the population dynamics.

## 4. Application

### 4.1 Problem setting

The problem considered in this section is a version of the benthic algae population dynamics considered in previous studies (Hamagami et al., 2024; Yoshioka and Hamagami, 2024)[38,72]. The problem setting is simple (**Figure 1**); a part of the riverbed of a river reach is assumed to be covered by nuisance benthic algal species, which are physically removed by sediment-laden (i.e., turbid) turbulent water flows while they grow by photosynthesis. Blooms of these benthic algae negatively affect aquatic food webs by extensively covering the riverbed; hence, investigations of their long-term influence on aquatic environments and ecosystems are important.

Yoshioka and Hamagami (2024) [38] reported that the physical removal of benthic algae can be modeled with long-memory decay, with the probability measure $p$ given by the gamma measure in **Remark 3**. Their model was a nonlinear SDE whose nonlinear drift contributes to generating long-memory decay, whereas the proposed model (5) is a superposition of linear SDEs without any nonlinearities and therefore serves as a more tractable mathematical framework. Moreover, previous studies did not discuss model uncertainty issues, which were approached in this study based on the proposed environmental index.

### 4.2 Parameter setting

We present a demonstrative application of the environmental index $V$ in (18) and optimal control in the sense of equilibrium control $\phi^*$ in (38), which is the worst-case model uncertainty maximizer in (18). In the remainder of this paper, we assume the most pessimistic case where the initial condition is $x_0 \equiv 1$, which, in the context of benthic algae population dynamics, implies that the riverbed is fully covered by algae.

We assume the gamma distributions discussed in **Remark 3** for $p$ and $\mu$ because these distributions are connected to the nonexponential discount as well as long memory, and gamma-type $p$ has been experimentally identified, as discussed in **Section 4.1**. Parameter $\alpha_\delta$ needs to be set larger than 1 if $R > 0$ according to **Proposition 1**. We assume the proportional case of $\eta(r) = cr$ with constant $c \in (0,1)$, where $c$ serves as an uncertainty aversion constant.

The parameters and coefficients are specified as follows. First, the probability measure $p$ was specified as follows: $\alpha_r = 0.295$ (-) and $\beta_r = 34.4$ (1/day) by hourly laboratory observations of algae removal in the observable area in an experimental flume (middle of **Figure 1**; see also **Section A3 in Appendix**). This is a truly long-memory case because $\alpha_r < 1$, and hence, the environmental index $V$ diverges if no discount exists. The growth rate $R$ is assumed to be $O(10^{-2})$ (1/day) based on Yoshioka and Hamagami (2024)[38]. We use 0.01 (1/day) or 0.02 (1/day). We have examined other values of $R$ but obtained qualitatively the same results (**Section A4 in the Appendix**). Concerning the discount factor, we set $\beta_\delta = 0.01$ (1/day) by assuming that $\beta_\delta << \beta_r$ because we want to consider cases where the discount is not small. The dimension of $r$ is (1/day).

Under these settings, we compute the environmental index $V$, which is the value function, and its maximizer $\phi^*$ for different values of $\alpha_\delta$, $R$, and $c$. The equation (36) is solved via a common fixed-point iteration. Because the coefficient $A$ in **Proposition 1** does not seem to be determined analytically, we compute it on the basis of the quantile discretization method specialized for gamma distributions along with importance sampling (Yoshioka and Yoshioka, 2024b)[74], where the sampling of discretization points $r_i$ and $\delta_j$ in **Section A2 in Appendix** is performed by using a quantile-based quantization method. We do not provide details of this discretization method but are available in Yoshioka and Yoshioka (2024b)[74]. We perform a convergence analysis in **Section A2**. We selected the discretization parameter values $N = M = 1{,}024$, which were sufficiently large according to convergence analysis.

### 4.3 Results and discussion

First, we investigate $\phi^*$ for different parameter settings. **Figure 2** shows the computed $\phi^*$, which depends on $R = 0.01$ and $\alpha_\delta = 1.25,\ 1.50,\ 2.00,\ 4.00$. Similarly, **Figure 3** shows the results of $R = 0.02$. Here, the limit $c \to 0$ corresponds to the non-robust case.

According to **Figures 2 and 3**, the worst-case model uncertainty $\phi^*$ decreases for the decay rate $r$ and decreases for the uncertainty aversion constant $c$. The dependence on $c$ is a reasonable property to be equipped with an environmental index because a more uncertainty-averse observer would evaluate the population dynamics in a more pessimistic way: in our context, a smaller $r\phi^*(r)$ and, hence, slower decay of the algae population. **Figures 2 and 3** quantify the extent to which the model should be distorted under the assumed degree of model uncertainty. The difference caused by the growth rate $R$ is not clearly visible at this stage but becomes more visible in the sequel.

Second, we discuss $V(1)$, the environmental index with the maximum initial population $x_0 \equiv 1$, which is truly the worst case among all the initial conditions. **Figure 4** shows the computed $V(1)$ for $R = 0.01$ and different values of $\alpha_\delta$. Similarly, **Figure 5** shows the computed $V(1)$ for $R = 0.02$. A comparison of **Figures 4 and 5** reveals that a larger growth rate yields $R$ larger $V(1)$ values parameterized through the uncertainty aversion parameter $c$, which is a desired property to be equipped with an environmental index, as it should be more pessimistic about model uncertainties. For both **Figures 4 and 5**, a smaller $\alpha_\delta$ value, that is, a smaller discount factor, results in a larger $V(1)$ value and slope, indicating that specifying a smaller discount factor yields a more sensitive environmental index. This is a desired property to be equipped with an environmental index, as it should be more pessimistic about model uncertainties.

We also discuss the efficient frontier shaped by the following two terms (43) and (44), which allows for a deeper analysis of our environmental index $V$. Here, we focus on the two terms contained in $V$, which represents the worst-case environmental impact

$$\mathrm{Env} = \int_0^{+\infty}\int_0^{+\infty} e^{-\delta s} X_s \mathrm{d}s\mu(\mathrm{d}\delta)\Big|_{\phi=\phi^*} \tag{43}$$

and the worst-case (discount-weighted) relative entropy

$$\mathrm{Ren} = \int_0^{+\infty}\int_0^{+\infty} e^{-\delta s}\left(\int_0^{+\infty} x_s(r)\big(\phi(r)\ln\phi(r) - \phi(r) + 1\big)p(\mathrm{d}r)\right)\mathrm{d}s\mu(\mathrm{d}\delta)\Big|_{\phi=\phi^*}, \tag{44}$$

where $V = \mathrm{Env} - \frac{1}{c}\mathrm{Ren}$ by definition. Both $\mathrm{Env}$ and $\mathrm{Ren}$ depend on $c$ through $\phi^*$; hence, we can obtain an efficient frontier as a curve parameterized by $c$ in a two-dimensional (2D) $(\mathrm{Ren}, \mathrm{Env})$ space.

**Figure 6** shows the computed $\mathrm{Ren}$-$\mathrm{Env}$ relationships for $R = 0.01$ and the different values of $\alpha_\delta$ parameterized by the uncertainty-aversion constant $c \in (0,1)$. Similarly, **Figure 7** shows the results of $R = 0.02$. For each fixed $R$ and $\alpha_\delta$, the $\mathrm{Ren}$-$\mathrm{Env}$ relationship yields concave and increasing plots, and both $\mathrm{Ren}$ and $\mathrm{Env}$ are increasing functions of $c$, which is in accordance with **Figures 2 and 3**, where the worst-case uncertainty $\phi^*$ decreases with respect to $c$ at each $r > 0$. Indeed, decreasing $\phi^*$ yields a larger $\phi^*\ln\phi^* - \phi^* + 1$ and a temporally slower decay of the population $x$. A finer inspection of the concavity of the $\mathrm{Ren}$-$\mathrm{Env}$ relationships implies that their increasing speed becomes constant as $\mathrm{Ren}$ increases. Hence, a nearly affine structure exists between $\mathrm{Ren}$ and $\mathrm{Env}$ for relatively large model uncertainties. These findings imply that the proposed environmental index, although infinite-dimensional, possesses a tractable geometric structure between the uncertainties and the population. Moreover, because **Figures 2 and 3** quantify how the model is distorted for each $c$, the results of **Figures 6 and 7** provide additional information to connect $c$ and $\mathrm{Ren}$ one-to-one at least computationally. Thus, our mathematical framework consistently links the assumed uncertainty size ( $\mathrm{Ren}$ ), proper overestimation of the disutility caused by the population ( $\mathrm{Env}$ ), and corresponding model distortion ( $\phi^{*}$ ). A comparison of **Figures 6 and 7** shows that a larger growth rate yields $R$ with the $\mathrm{Ren}$-$\mathrm{Env}$ relationship shifting vertically upward in the 2D $(\mathrm{Ren}, \mathrm{Env})$ space for the environmental index $V$, suggesting that the efficient frontier shaped by the two indices $\mathrm{Ren}$ and $\mathrm{Env}$ can serve as an alternative environmental index to $V$. An analysis that combines these quantities potentially provides finer information than studying them individually, as demonstrated in this study.

Our approach assumes a long-memory process as a target system while also covering classical exponential-memory processes because the latter can be realized by assuming a probability measure $\mu$ concentrated at a positive value. The population dynamics in an engineering problem are nonlinear (e.g., Maurya and Misra, 2025; Sajid et al., 2024)[75,76]; thus, the proposed approach does not directly apply to them; however, if the problem of interest is the population dynamics around an equilibrium point, then the system can be linearized such that the proposed approach performs well. Such cases can be accommodated by suitably tuning $\mu$ and $R$ and then computing the environmental index numerically, as demonstrated in this section. In practice, reducing model uncertainties would be possible via frequent and automatic observation schemes (Beal and Schaeffer, 2026)[77], but their operational costs may have to be accounted for in the model, which was not considered in the proposed framework. Nevertheless, its control-theoretical nature could be helpful in developing such an advanced model.

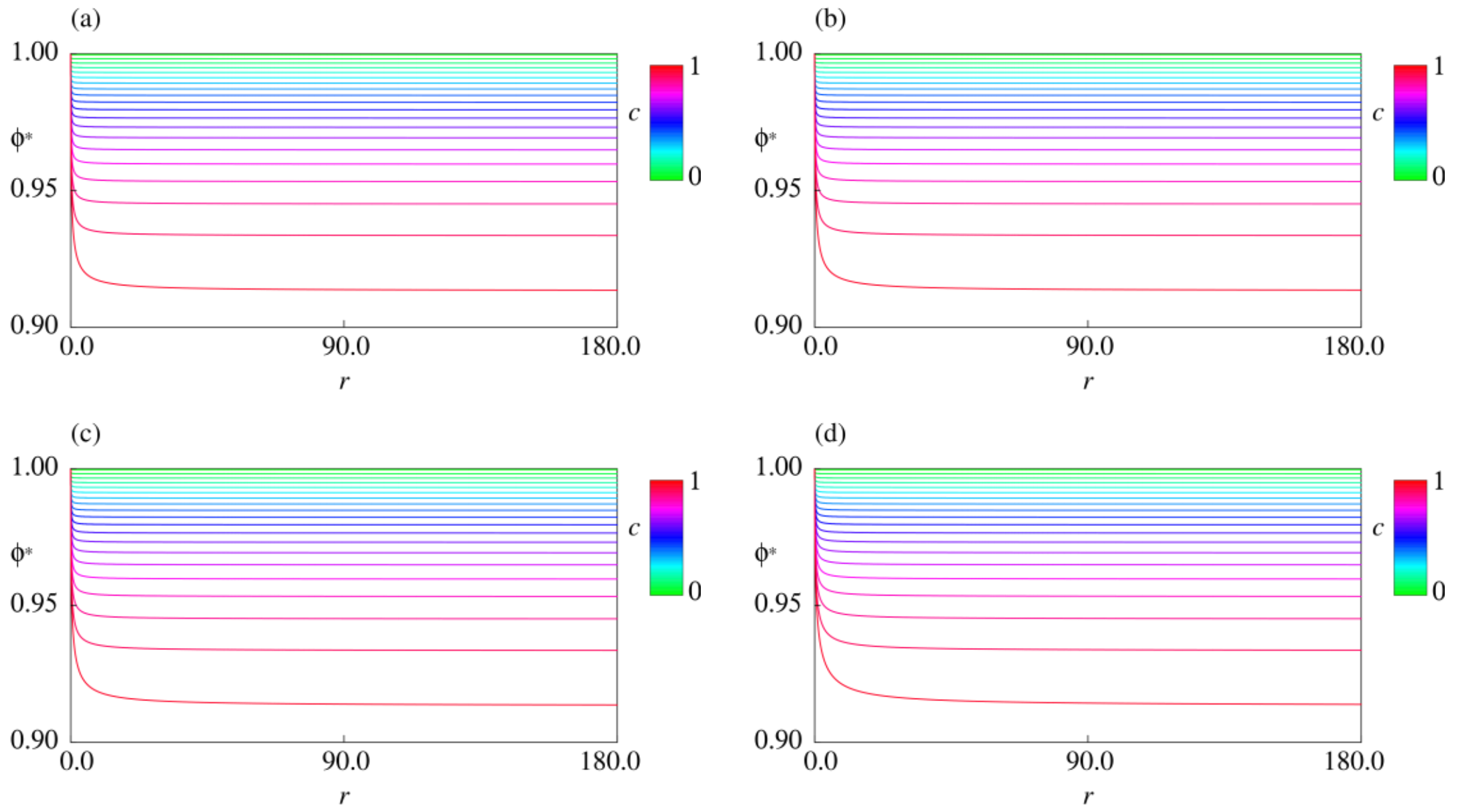

**Figure 2.** Computed $\phi^*$ for $R = 0.01$ and different values of $\alpha_\delta$: (a) 1.25, (b) 1.50, (c) 2.00, and (d) 4.00. Each of the figure panels may be apparently not critically different, but their difference becomes apparent in the corresponding environmental indices; similar observations apply to **Figure 3**.

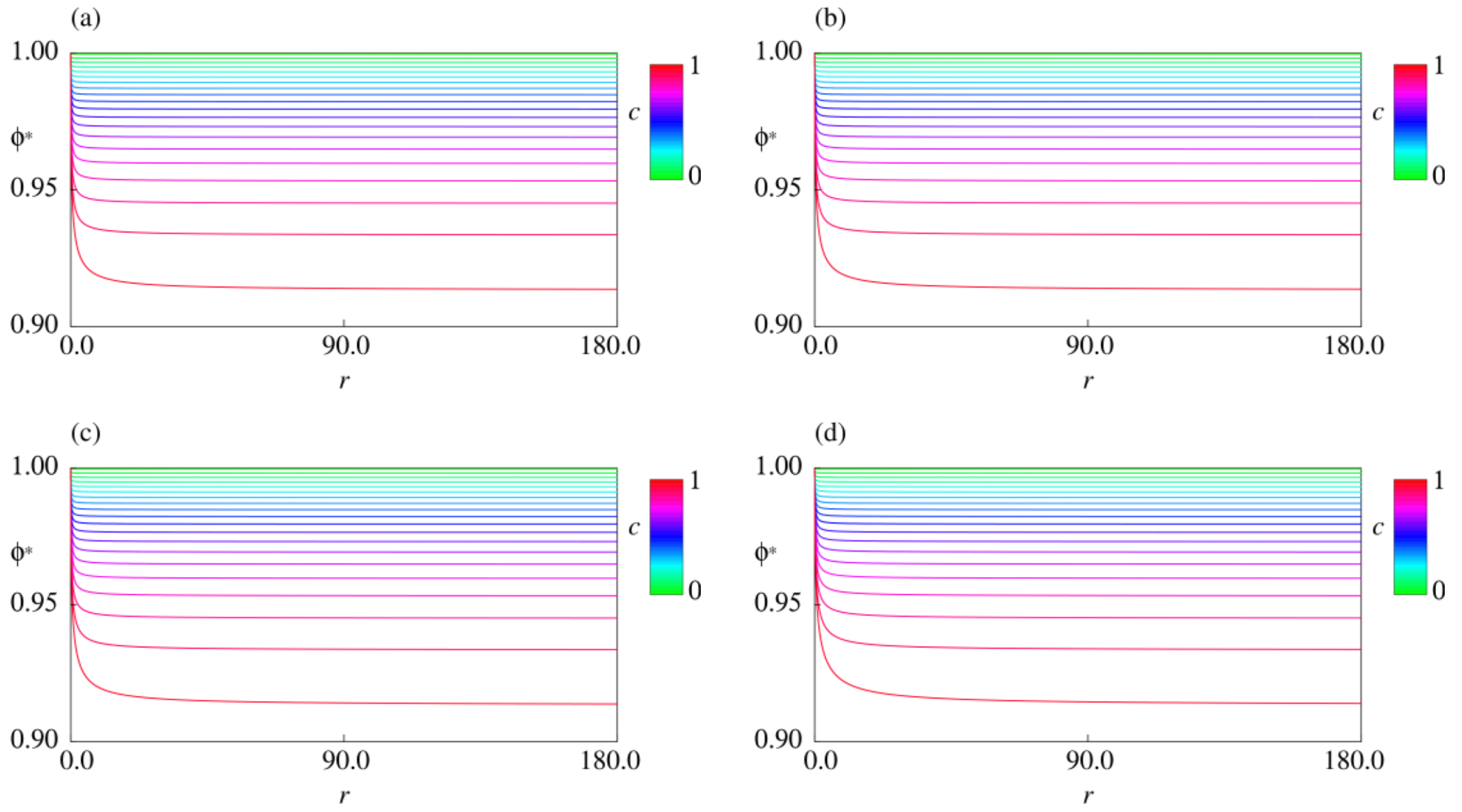

**Figure 3.** Computed $\phi^*$ for $R = 0.02$ and different values of $\alpha_\delta$: (a) 1.25, (b) 1.50, (c) 2.00, and (d) 4.00.

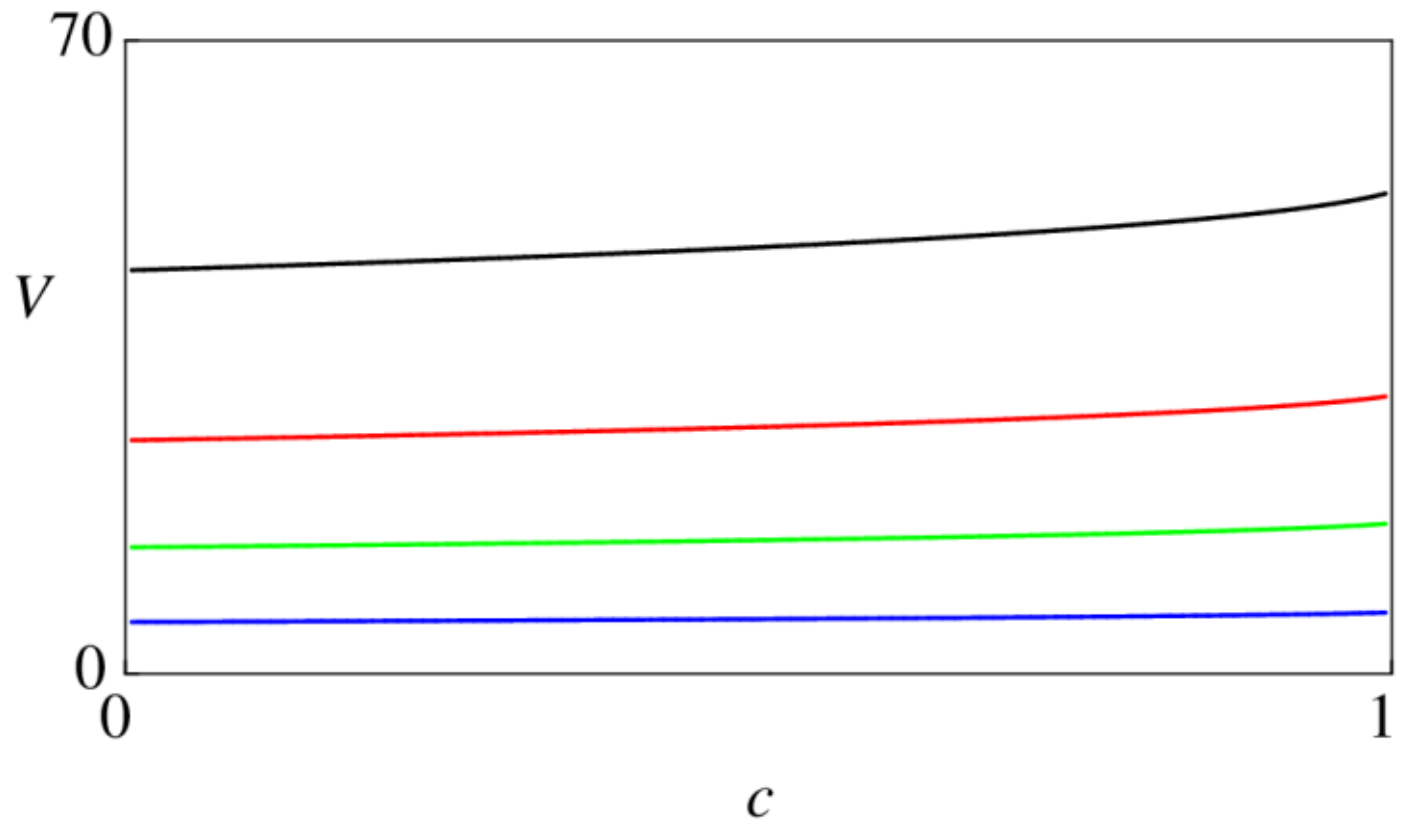

**Figure 4.** Computed $V = V(1)$ for $R = 0.01$ and different values of $\alpha_\delta$: 1.25 (black), 1.50 (red), 2.00 (green), and 4.00 (blue).

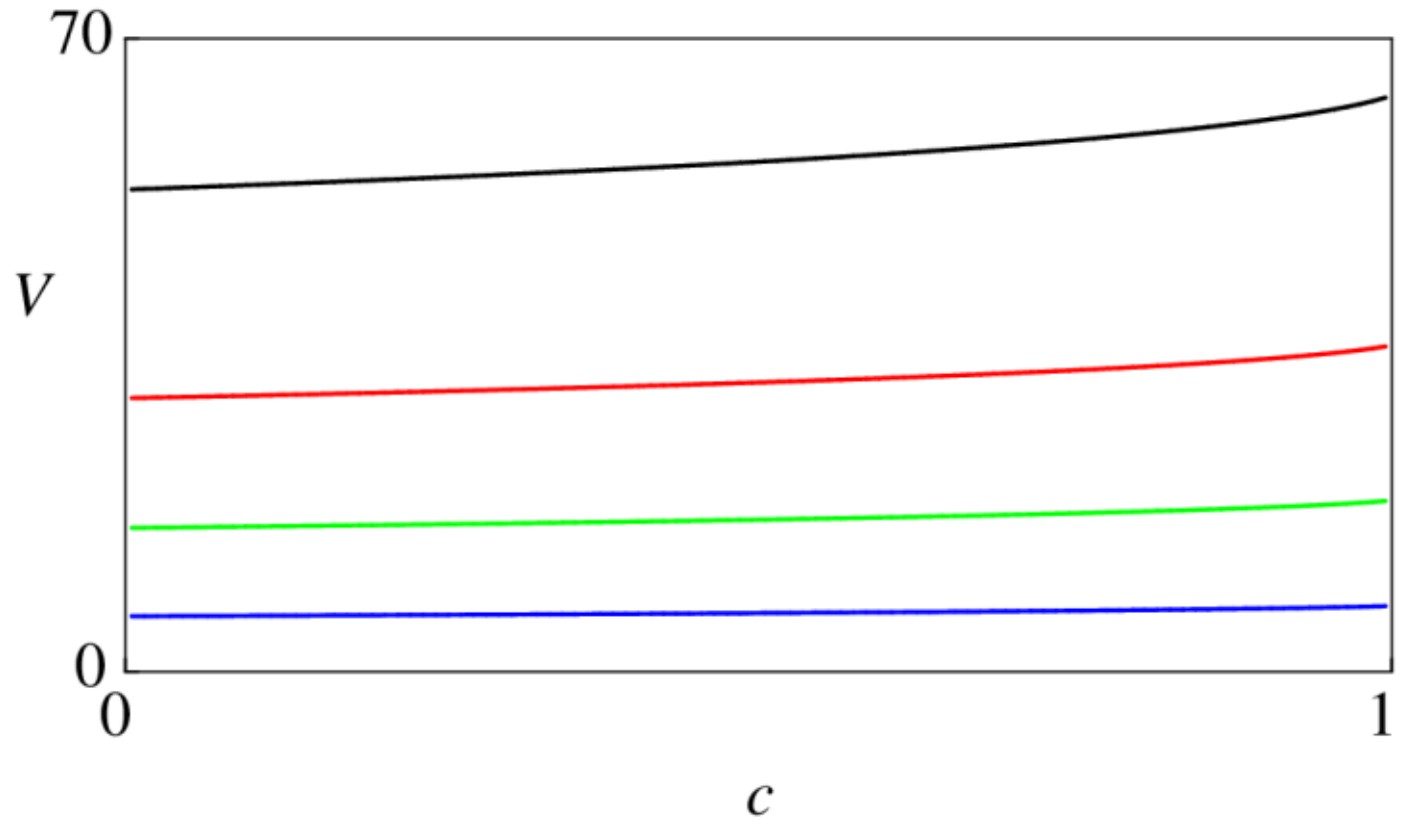

**Figure 5.** Computed $V = V(1)$ for $R = 0.02$ and different values of $\alpha_\delta$: 1.25 (black), 1.50 (red), 2.00 (green), and 4.00 (blue).

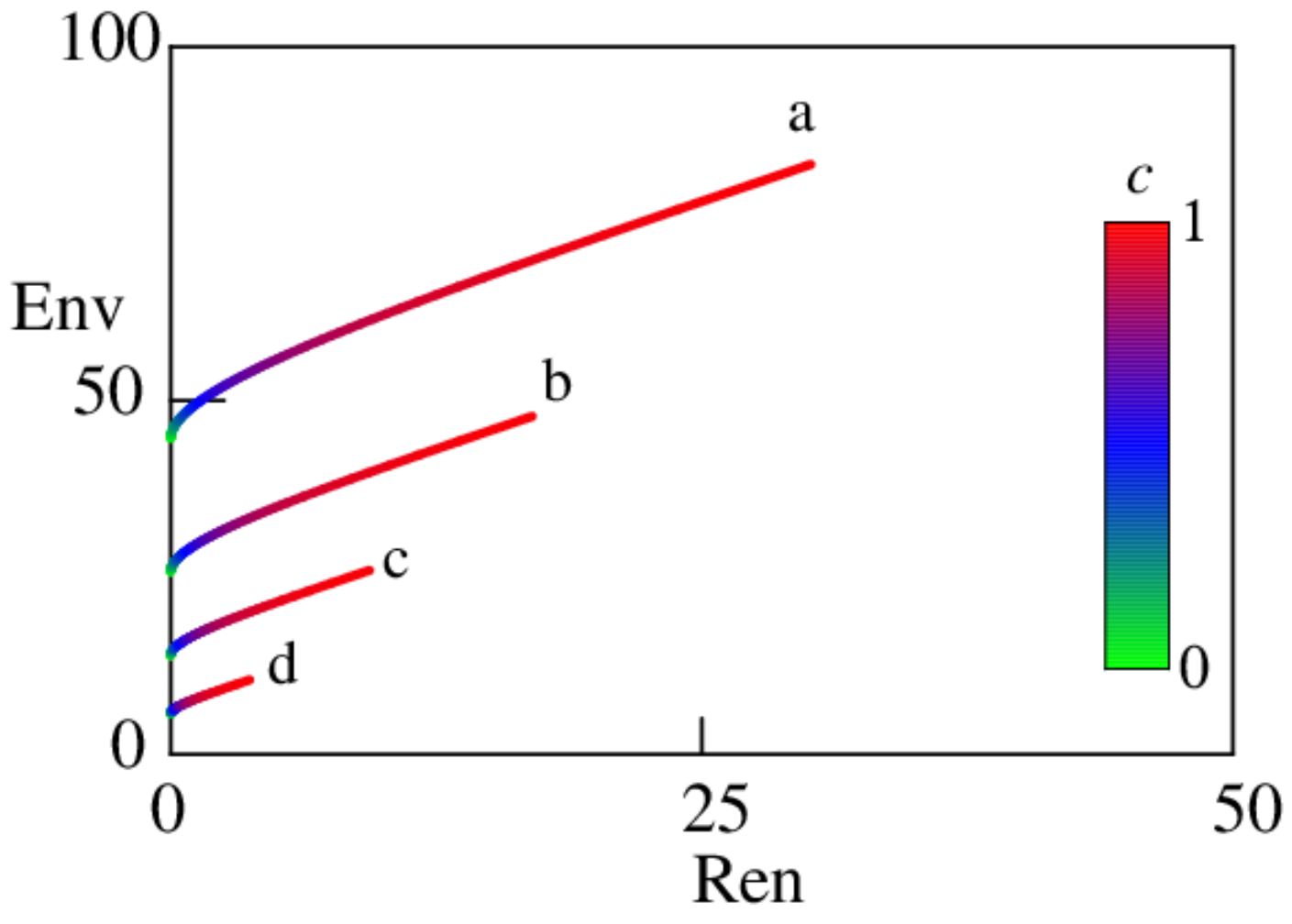

**Figure 6.** Computed Ren-Env relationship for $R = 0.01$ and different values of $\alpha_\delta$ : (a) 1.25, (b) 1.50, (c) 2.00, and (d) 4.00.

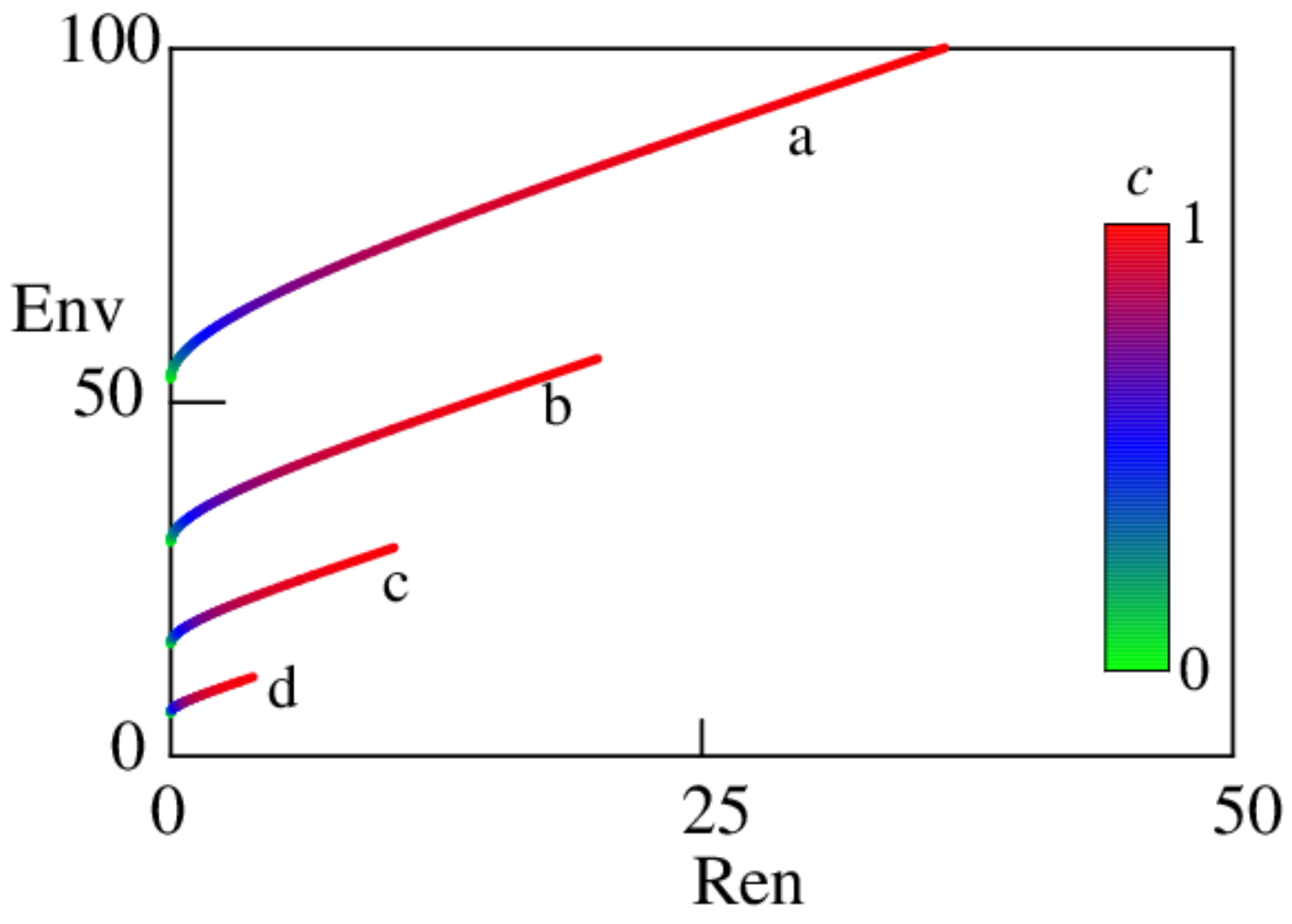

**Figure 7.** Computed Ren-Env relationship for $R = 0.02$ and different values of $\alpha_\delta$ : (a) 1.25, (b) 1.50, (c) 2.00, and (d) 4.00.

## 5. Conclusion

We proposed a tractable mathematical model for assessing the cumulative environmental impacts caused by long-memory processes subject to model uncertainties. A semi-explicit solution to the extended HJB system associated with the time-inconsistent control problem yielded an environmental index whose existence is guaranteed under certain conditions of a nonexponential discount. A computable finite-dimensional version of the solution was presented and applied to a demonstrative example focusing on benthic algae management in river environments. The proposed mathematical approach paves the way for a new environmental impact assessment for persistent environmental factors on the basis of long-memory processes combined with time-inconsistent control theory. Moreover, we demonstrated that even simple population dynamics could lead to a nontrivial environmental index, which is analytically tractable.

We did not address the assessment of time-periodic environmental factors, which would serve as a more realistic mathematical model than that used in this study, because environmental variables are influenced by seasonal and/or diurnal perturbations. Another research topic that we did not study was the impact assessment of nonlinear long-memory processes, with which the corresponding extended HJB system does not have any semi-explicit solutions. Our solution was smooth, whereas in such a case, a proper solution would not necessarily be continuously differentiable everywhere, and the notion of viscosity solutions can become an effective tool. From an application perspective, collecting data on (nonlinear) long-memory processes is an important task that will be continued by the authors, which is based on both experimental and field data collection. Finally, model uncertainties can be evaluated via a generalized divergence that truly extends the relative entropy (Davis et al., 2025)[26]. We expect that the proposed framework can be applied to this case if the admissibility condition for the guessed solution is modified, which will be addressed in future research.

## Appendix

### A1. Proofs

***Proof of Proposition 1***

The proof is divided into several sections. For any $f \in E$, we obtain the following variational derivatives:

$$\Phi_f(r, f) = A(r) \tag{45}$$

and

$$\Psi_{f,t}(r, f) = \int_0^{+\infty} C(r, \delta) e^{-\delta t} \mu(\mathrm{d}\delta). \tag{46}$$

**Step 1: Equation of** $A$

First, we show that (36) admits a unique solution in the class of $p$ a.e. positive functions. We consider (36) for an arbitrary $r > 0$ in the remainder of **Step 1**. We can rewrite (36) as

$$A = F(A), \quad A \geq 0 \tag{47}$$

with

$$F(A) = I(A)\left(1 - r\frac{1 - (1 + A\eta(r))e^{-\eta(r)A}}{\eta(r)}\right) \tag{48}$$

and

$$I(A) = \int_0^{+\infty} \frac{1}{\delta + re^{-\eta(r)A} + R} \mu(\mathrm{d}\delta) > 0. \tag{49}$$

We have

$$re^{-\eta(r)A} I(A) = \int_0^{+\infty} \frac{re^{-\eta(r)A}}{\delta + re^{-\eta(r)A} + R} \mu(\mathrm{d}\delta) \in (0,1) \tag{50}$$

and

$$\frac{\mathrm{d}}{\mathrm{d}A} I(A) = \int_0^{+\infty} \frac{\eta(r) re^{-\eta(r)A}}{\left(\delta + re^{-\eta(r)A} + R\right)^2} \mu(\mathrm{d}\delta) < \eta(r) I(A), \quad A \geq 0. \tag{51}$$

The elemental calculations show

$$0 - F(0) = 0 - I(0) < 0 \tag{52}$$

and

$$
\begin{aligned}
\bar{A}(r)-F(\bar{A}(r)) &= \bar{A}(r)-I(\bar{A}(r))\left(1-r\frac{1-(1+\bar{A}(r)\eta(r))e^{-\eta(r)\bar{A}(r)}}{\eta(r)}\right) \\
&= \bar{A}(r)-I(\bar{A}(r))\left(1-r\frac{1-e^{-\eta(r)\bar{A}(r)}}{\eta(r)}+r\bar{A}(r)e^{-\eta(r)\bar{A}(r)}\right) \\
&= \bar{A}(r)-I(\bar{A}(r))r\bar{A}(r)e^{-\eta(r)\bar{A}(r)} \\
&= \bar{A}(r)\left(1-re^{-\eta(r)\bar{A}(r)}I(\bar{A}(r))\right) \\
&> 0
\end{aligned} \tag{53}
$$

from (51) and the definition of $\bar{A}(r)$. For $0<A\le\bar{A}(r)$, we obtain

$$
\begin{aligned}
\frac{\mathrm{d}}{\mathrm{d}A}F(A) &= \frac{\mathrm{d}}{\mathrm{d}A}I(A)\left(1-r\frac{1-(1+A\eta(r))e^{-\eta(r)A}}{\eta(r)}\right)+I(A)\frac{\mathrm{d}}{\mathrm{d}A}\left(1-r\frac{1-(1+A\eta(r))e^{-\eta(r)A}}{\eta(r)}\right) \\
&= \frac{\mathrm{d}}{\mathrm{d}A}I(A)\left(1-r\frac{1-(1+A\eta(r))e^{-\eta(r)A}}{\eta(r)}\right)-I(A)r\eta(r)Ae^{-\eta(r)A} \\
&< \eta(r)I(A)\left(1-r\frac{1-(1+A\eta(r))e^{-\eta(r)A}}{\eta(r)}\right)-I(A)r\eta(r)Ae^{-\eta(r)A} \\
&= I(A)\left(\eta(r)-r+r(1+A\eta(r))e^{-\eta(r)A}\right)-I(A)r\eta(r)Ae^{-\eta(r)A} \\
&= I(A)\left(\eta(r)-r+re^{-\eta(r)A}\right) \\
&< I(A)\left(r-r+re^{-\eta(r)A}\right) \\
&= I(A)re^{-\eta(r)A} \\
&< 1,
\end{aligned} \tag{54}
$$

where we used (50) and

$$
1-r\frac{1-(1+A\eta(r))e^{-\eta(r)A}}{\eta(r)}=1-r\frac{1-e^{-\eta(r)A}}{\eta(r)}+Are^{-\eta(r)A}\ge 1-r\frac{1-e^{-\eta(r)\bar{A}(r)}}{\eta(r)}+Are^{-\eta(r)A}>0 \tag{55}
$$

for $0\le A\le\bar{A}(r)$.

Consequently, the function $A-F(A)$ is increasing for $0\le A\le\bar{A}(r)$ and is continuous (this fact directly follows from the forms of (48) and (49)). This observation, along with (52) and (53), shows that (36) admits a unique positive solution $p$ a.e., satisfying the bound (37). Finally, as a byproduct, we have the following bound (see also the beginning of **Step 2**):

$$
0\le\frac{r}{\eta(r)}f(r)\left(\phi^*(r)\ln\phi^*(r)-\phi^*(r)+1\right)=\frac{r}{\eta(r)}\left(1-(1+A(r)\eta(r))e^{-\eta(r)A(r)}\right)\le 1\,. \tag{56}
$$

***Step 2: Guessed solution: formal check***

Substituting (27) into (25) yields

$$
\phi^*(\cdot)=\exp(-\eta(\cdot)A(\cdot)), \tag{57}
$$

which is bounded between 0 and 1 because the quantity inside "exp" is not positive. For any $f \in E$,

$$\sup_{\phi>0}\left\{\begin{array}{l}-\int_0^{+\infty} r\phi(r)f(r)\Phi_f(r,f)p(\mathrm{d}r)\\ -\int_0^{+\infty}\frac{r}{\eta(r)}f(r)\left(\phi(r)\ln\phi(r)-\phi(r)+1\right)p(\mathrm{d}r)\end{array}\right\} = -\int_0^{+\infty} f(r)\frac{r}{\eta(r)}\left(1-e^{-\eta(r)A(r)}\right)p(\mathrm{d}r). \quad (58)$$

Substituting (27) and (28) into (23) yields

$$\begin{aligned}&-\int_0^{+\infty} f(r)\frac{r}{\eta(r)}\left(1-e^{-\eta(r)A(r)}\right)p(\mathrm{d}r)\\ &+\int_0^{+\infty} R\left(1-f(r)\right)A(r)p(\mathrm{d}r)+\int_0^{+\infty} f(r)p(\mathrm{d}r)+\int_0^{+\infty}\left\{\int_0^{+\infty} C(r,\delta)f(r)+D(\delta)\right\}\mu(\mathrm{d}\delta)=0\end{aligned}. \quad (59)$$

Similarly, substituting (27) and (28) into (23) yields

$$\begin{aligned}&-\int_0^{+\infty}\delta\left\{\int_0^{+\infty} C(r,\delta)f(r)+D(\delta)\right\}e^{-\delta t}\mu(\mathrm{d}\delta)\\ &-\int_0^{+\infty} r\phi^*(r)f(r)\int_0^{+\infty} C(r,\delta)e^{-\delta t}\mu(\mathrm{d}\delta)p(\mathrm{d}r)+\int_0^{+\infty} R\left(1-f(r)\right)\int_0^{+\infty} C(r,\delta)e^{-\delta t}\mu(\mathrm{d}\delta)p(\mathrm{d}r)\\ &+\rho_t\left\{\int_0^{+\infty} f(r)p(\mathrm{d}r)-\int_0^{+\infty}\frac{r}{\eta(r)}f(r)\left(\phi^*(r)\ln\phi^*(r)-\phi^*(r)+1\right)p(\mathrm{d}r)\right\}=0.\end{aligned} \quad (60)$$

From (59), because $f$ is arbitrary, we obtain

$$-\frac{r}{\eta(r)}\left(1-e^{-\eta(r)A(r)}\right)-RA(r)+1+\int_0^{+\infty} C(r,\delta)\mu(\mathrm{d}\delta)=0,\ \ p \text{ a.e.} \quad (61)$$

and

$$R\int_0^{+\infty} A(r)p(\mathrm{d}r)+\int_0^{+\infty} D(\delta)\mu(\mathrm{d}\delta)=0. \quad (62)$$

From (60), we subsequently obtain (we used (21))

$$-\delta C(r,\delta)-re^{-\eta(r)A(r)}C(r,\delta)-RC(r,\delta)-\delta\left(1-\frac{r}{\eta(r)}\left(1-\left(1+\eta(r)A(r)\right)e^{-\eta(r)A(r)}\right)\right)=0 \quad p\times\mu \text{ a.e.} \quad (63)$$

and

$$-\delta D(\delta)+R\int_0^{+\infty} C(r,\delta)p(\mathrm{d}r)=0 \quad \mu \text{ a.e.} \quad (64)$$

Rearranging (61), (63), and (64) yields (36), (34), and (35), respectively.

The remaining coefficient $B$ is not found from (61)–(64). From the guessed solution, we have

$$\begin{aligned}\Phi(x_0)&=\int_0^{+\infty} A(r)x_0(r)p(\mathrm{d}r)+B\\ &=\int_0^{+\infty}\int_0^{+\infty} e^{-\delta s}\int_0^{+\infty} x_s(r)p(\mathrm{d}r)\mathrm{d}s\mu(\mathrm{d}\delta)\\ &-\int_0^{+\infty}\int_0^{+\infty} e^{-\delta s}\left(\int_0^{+\infty} x_s(r)\left(1-\frac{r}{\eta}\left(\phi^*(r)\ln\phi^*(r)-\phi^*(r)+1\right)\right)p(\mathrm{d}r)\right)\mathrm{d}s\mu(\mathrm{d}\delta)\\ &=\int_0^{+\infty}\int_0^{+\infty}\left\{\begin{array}{l}\left(1-\frac{r}{\eta(r)}\left(\phi^*(r)\ln\phi^*(r)-\phi^*(r)+1\right)\right)\\ \times\left(\frac{R}{\delta\left(R+r\phi^*(r)\right)}+\frac{1}{\delta+R+r\phi^*(r)}\left(x_0(r)-\frac{R}{R+r\phi^*(r)}\right)\right)p(\mathrm{d}r)\end{array}\right\}\mu(\mathrm{d}\delta),\end{aligned} \quad (65)$$

from which we obtain

$$
\begin{aligned}
B &= \int_0^{+\infty}\int_0^{+\infty}\left\{\begin{aligned}&\left(1-\frac{r}{\eta(r)}\left(\phi^*(r)\ln\phi^*(r)-\phi^*(r)+1\right)\right)\\ &\times\left(\frac{R}{\delta\left(R+r\phi^*(r)\right)}+\frac{1}{\delta+R+r\phi^*(r)}\left(-\frac{R}{R+r\phi^*(r)}\right)\right)p(\mathrm{d}r)\end{aligned}\right\}\mu(\mathrm{d}\delta)\\
&= R\int_0^{+\infty}\int_0^{+\infty}\left\{\left(1-\frac{r}{\eta(r)}\left(\phi^*(r)\ln\phi^*(r)-\phi^*(r)+1\right)\right)\frac{1}{\delta\left(\delta+R+r\phi^*(r)\right)}p(\mathrm{d}r)\right\}\mu(\mathrm{d}\delta) \qquad (66)\\
&= R\int_0^{+\infty}\int_0^{+\infty}\left\{\left(1-\frac{r}{\eta(r)}\left(1-\left(1+\eta(r)A(r)\right)e^{-\eta(r)A(r)}\right)\right)\frac{1}{\delta\left(\delta+R+r\phi^*(r)\right)}p(\mathrm{d}r)\right\}\mu(\mathrm{d}\delta).
\end{aligned}
$$

This gives (33).

We check the terminal condition (24). From (28) and (35), for any $f\in E$, we have

$$
\begin{aligned}
&\Psi_t(f)\\
&=\int_0^{+\infty}\left\{\int_0^{+\infty}C(r,\delta)f(r)p(\mathrm{d}r)+\frac{R}{\delta}\int_0^{+\infty}C(r,\delta)p(\mathrm{d}r)\right\}e^{-\delta t}\mu(\mathrm{d}\delta)\\
&=\int_0^{+\infty}\int_0^{+\infty}\left(\frac{R}{\delta}+f(r)\right)C(r,\delta)p(\mathrm{d}r)e^{-\delta t}\mu(\mathrm{d}\delta) \qquad (67)\\
&=-\int_0^{+\infty}\int_0^{+\infty}\left(\frac{R}{\delta}+f(r)\right)\frac{\delta}{\delta+re^{-\eta(r)A(r)}+R}\left(1-\frac{r}{\eta(r)}\left(1-\left(1+\eta(r)A(r)\right)e^{-\eta(r)A(r)}\right)\right)p(\mathrm{d}r)e^{-\delta t}\mu(\mathrm{d}\delta).
\end{aligned}
$$

From (55), we obtain

$$
\begin{aligned}
&\left|\Psi_t(f)\right|\\
&\le\int_0^{+\infty}\int_0^{+\infty}\left(\frac{R}{\delta}+\left|f(r)\right|\right)\frac{\delta}{\delta+re^{-\eta(r)A(r)}+R}\left(1-\frac{r}{\eta(r)}\left(1-\left(1+\eta(r)A(r)\right)e^{-\eta(r)A(r)}\right)\right)p(\mathrm{d}r)e^{-\delta t}\mu(\mathrm{d}\delta), \qquad (68)
\end{aligned}
$$

where we used (37). To prove (24), it suffices to show the integrability condition because $f$ is bounded:

$$
\int_0^{+\infty}\int_0^{+\infty}\left(\frac{R}{\delta}+1\right)\frac{\delta}{\delta+r-\eta(r)+R}\left(1-\frac{r}{\eta(r)}\left(1-\left(1+\eta(r)A(r)\right)e^{-\eta(r)A(r)}\right)\right)p(\mathrm{d}r)\mu(\mathrm{d}\delta)<+\infty. \qquad (69)
$$

Therefore, we analyze (69). The inequality (56) implies that (69) is satisfied if

$$
\int_0^{+\infty}\int_0^{+\infty}\left(\frac{R}{\delta}+1\right)\frac{\delta}{\delta+r-\eta(r)+R}p(\mathrm{d}r)\mu(\mathrm{d}\delta)<+\infty. \qquad (70)
$$

We address the cases $R=0$ and $R>0$ separately. For $R=0$, (70) becomes

$$
\int_0^{+\infty}\int_0^{+\infty}\frac{\delta}{\delta+r-\eta(r)}p(\mathrm{d}r)\mu(\mathrm{d}\delta)<+\infty. \qquad (71)
$$

This condition is automatically satisfied because the integrand of (71) is between 0 and 1 since $r-\eta(r)>0$. For $R>0$, (70) is satisfied if

$$
\int_0^{+\infty}\int_0^{+\infty}\frac{1}{\delta+r-\eta(r)+R}p(\mathrm{d}r)\mu(\mathrm{d}\delta)<+\infty, \qquad (72)
$$

but this inequality holds true because $R>0$.

**Step 2** is completed by proving that (62) is redundant and follows from (61), (63), and (64). If this is completed, it follows that the guessed solution satisfies the extended HJB. We write

$$H(r)=1-\frac{r}{\eta(r)}\left(1-\left(1+\eta(r)A(r)\right)e^{-\eta(r)A(r)}\right). \tag{73}$$

From (61), we have

$$\left(re^{-\eta(r)A(r)}+R\right)A(r)=H(r)+\int_0^{+\infty}C(r,\delta)\mu(\mathrm{d}\delta). \tag{74}$$

Then, we obtain the desired equality via direct substitution:

$$\begin{aligned}
&\frac{1}{R}\left\{R\int_0^{+\infty}A(r)p(\mathrm{d}r)+\int_0^{+\infty}D(\delta)\mu(\mathrm{d}\delta)\right\}\\
&=\int_0^{+\infty}\frac{1}{re^{-\eta(r)A(r)}+R}H(r)p(\mathrm{d}r)+\int_0^{+\infty}\int_0^{+\infty}\left(\frac{1}{re^{-\eta(r)A(r)}+R}+\frac{1}{\delta}\right)C(r,\delta)\mu(\mathrm{d}\delta)p(\mathrm{d}r)\\
&=\int_0^{+\infty}\frac{1}{re^{-\eta(r)A(r)}+R}H(r)p(\mathrm{d}r)+\int_0^{+\infty}\int_0^{+\infty}\left(\frac{\delta+re^{-\eta(r)A(r)}+R}{\delta\left(re^{-\eta(r)A(r)}+R\right)}\right)C(r,\delta)\mu(\mathrm{d}\delta)p(\mathrm{d}r)\\
&=\int_0^{+\infty}\frac{1}{re^{-\eta(r)A(r)}+R}H(r)p(\mathrm{d}r)\\
&-\int_0^{+\infty}\int_0^{+\infty}\left(\frac{\delta+re^{-\eta(r)A(r)}+R}{\delta\left(re^{-\eta(r)A(r)}+R\right)}\right)\frac{\delta}{\delta+re^{-\eta(r)A(r)}+R}H(r)\mu(\mathrm{d}\delta)p(\mathrm{d}r)\\
&=\int_0^{+\infty}\frac{1}{re^{-\eta(r)A(r)}+R}H(r)p(\mathrm{d}r)-\int_0^{+\infty}\int_0^{+\infty}\frac{1}{re^{-\eta(r)A(r)}+R}H(r)\mu(\mathrm{d}\delta)p(\mathrm{d}r)\\
&=\int_0^{+\infty}\frac{1}{re^{-\eta(r)A(r)}+R}H(r)p(\mathrm{d}r)-\int_0^{+\infty}\frac{1}{re^{-\eta(r)A(r)}+R}H(r)p(\mathrm{d}r)\\
&=0.
\end{aligned} \tag{75}$$

**Step 3: Guessed solution: boundedness**

We check the boundedness condition (29) for each coefficient. It suffices to check this condition for $A$ and $B$ because that for $C$ and $D$ has already been verified in **Step 2**.

***Coefficient** A*

We must show

$$\int_0^{+\infty}A(r)p(\mathrm{d}r)<+\infty, \tag{76}$$

but this follows from the following strict boundedness of $A$:

$$
\begin{aligned}
A(r) &= I(A)\left(1 - r\underbrace{\frac{1-\left(1+A(r)\eta(r)\right)e^{-\eta(r)A(r)}}{\eta(r)}}_{\text{nonnegative}}\right) \\
&\le I(A) \\
&\le \int_0^{+\infty} \frac{1}{\delta+R}\mu(\mathrm{d}\delta) \\
&< +\infty.
\end{aligned} \tag{77}
$$

***Coefficient*** $B$

We must show that $|B| < +\infty$, which is satisfied if

$$
\int_0^{+\infty}\int_0^{+\infty} \frac{R}{\delta\left(\delta + R + re^{-\eta(r)A(r)}\right)}\mu(\mathrm{d}\delta)p(\mathrm{d}r) < +\infty. \tag{78}
$$

If $R=0$, this condition is trivially satisfied. If $R>0$, then (78) follows from (32) and

$$
\frac{R}{\delta\left(\delta + R + re^{-\eta(r)A(r)}\right)} < \frac{1}{\delta}. \tag{79}
$$

**Step 4: Admissibility conditions**

We check the admissibility conditions (16) and (17). First, for (16), using the control $\phi = \phi^*$ that is time independent, we have

$$
\begin{aligned}
&\int_0^{+\infty}\int_0^{+\infty} e^{-\delta s}X_s \mathrm{d}s\mu(\mathrm{d}\delta) \\
&= \int_0^{+\infty}\int_0^{+\infty} e^{-\delta s}\int_0^{+\infty}\left(\frac{R}{R+r\phi^*(r)} + \left(x_0(r) - \frac{R}{R+r\phi^*(r)}\right)e^{-\left(R+r\phi^*(r)\right)s}\right)p(\mathrm{d}r)\mathrm{d}s\mu(\mathrm{d}\delta) \\
&= \int_0^{+\infty}\int_0^{+\infty}\left(\frac{R}{\delta\left(R+r\phi^*(r)\right)} + \frac{1}{\delta+R+r\phi^*(r)}\left(x_0(r) - \frac{R}{R+r\phi^*(r)}\right)\right)p(\mathrm{d}r)\mu(\mathrm{d}\delta) \\
&\le \int_0^{+\infty}\int_0^{+\infty}\left(\frac{R}{\delta\left(R+r\phi^*(r)\right)} + \frac{1}{\delta+R+r\phi^*(r)}\left(1 - \frac{R}{R+r\phi^*(r)}\right)\right)p(\mathrm{d}r)\mu(\mathrm{d}\delta) \\
&= \int_0^{+\infty}\int_0^{+\infty}\left(\frac{R}{\delta\left(R+r\phi^*(r)\right)} + \frac{1}{\delta+R+r\phi^*(r)}\frac{r\phi^*(r)}{R+r\phi^*(r)}\right)p(\mathrm{d}r)\mu(\mathrm{d}\delta) \\
&\le \int_0^{+\infty}\int_0^{+\infty}\left(\frac{R}{\delta\left(R+r\phi^*(r)\right)} + \frac{1}{\delta+R+r\phi^*(r)}\right)p(\mathrm{d}r)\mu(\mathrm{d}\delta).
\end{aligned} \tag{80}
$$

If $R>0$, then

$$
\int_0^{+\infty}\int_0^{+\infty} e^{-\delta s}X_s \mathrm{d}s\mu(\mathrm{d}\delta) \le \int_0^{+\infty}\int_0^{+\infty}\left(\frac{1}{\delta}+\frac{1}{R}\right)p(\mathrm{d}r)\mu(\mathrm{d}\delta), \tag{81}
$$

and hence we obtain (16) ad for (78). If $R=0$, then

$$\begin{aligned}\int_0^{+\infty}\int_0^{+\infty} e^{-\delta s} X_s \mathrm{d}s\mu(\mathrm{d}\delta) &\le \int_0^{+\infty}\int_0^{+\infty}\left(\frac{1}{\delta + r\phi^*(r)}\right)p(\mathrm{d}r)\mu(\mathrm{d}\delta)\\ &\le \int_0^{+\infty}\int_0^{+\infty}\left(\frac{1}{\delta + r - \eta(r)}\right)p(\mathrm{d}r)\mu(\mathrm{d}\delta)\\ &\le \int_0^{+\infty}\int_0^{+\infty}\left(\frac{1}{\delta + (1-\bar{c})r}\right)p(\mathrm{d}r)\mu(\mathrm{d}\delta)\\ &\le \frac{1}{1-\bar{c}}\int_0^{+\infty}\int_0^{+\infty}\frac{1}{\delta + r}p(\mathrm{d}r)\mu(\mathrm{d}\delta),\end{aligned} \tag{82}$$

and hence we obtain (17) from (31).

Second, for (17), with $\phi = \phi^*$, we have

$$\begin{aligned}&\int_0^{+\infty}\int_0^{+\infty} e^{-\delta s}\left(\int_0^{+\infty} x_s(r)\left(1-\frac{r}{\eta}\left(\phi^*(r)\ln\phi^*(r)-\phi^*(r)+1\right)\right)p(\mathrm{d}r)\right)\mathrm{d}s\mu(\mathrm{d}\delta)\\ &= \int_0^{+\infty}\int_0^{+\infty}\left(1-\frac{r}{\eta}\left(\phi^*(r)\ln\phi^*(r)-\phi^*(r)+1\right)\right)\left(\begin{array}{c}\dfrac{R}{\delta\left(R+r\phi^*(r)\right)}\\ +\dfrac{1}{\delta+R+r\phi^*(r)}\left(x_0(r)-\dfrac{R}{R+r\phi^*(r)}\right)\end{array}\right)p(\mathrm{d}r)\mu(\mathrm{d}\delta)\\ &\le \int_0^{+\infty}\int_0^{+\infty}(1-0)\left(\frac{R}{\delta\left(R+r\phi^*(r)\right)}+\frac{1}{\delta+R+r\phi^*(r)}\left(x_0(r)-\frac{R}{R+r\phi^*(r)}\right)\right)p(\mathrm{d}r)\mu(\mathrm{d}\delta)\\ &\le \int_0^{+\infty}\int_0^{+\infty}\left(\frac{R}{\delta\left(R+r\phi^*(r)\right)}+\frac{1}{\delta+R+r\phi^*(r)}\frac{r\phi^*(r)}{R+r\phi^*(r)}\right)p(\mathrm{d}r)\mu(\mathrm{d}\delta),\end{aligned} \tag{83}$$

where we used the fact that $\phi^*$ is bounded between 0 and 1 and that the function $y\ln y - y + 1$ ($0 \le y \le 1$) is nonnegative and maximized at 0 with a maximum value of 1. Then, we obtain (17) when $\phi = \phi^*$, as discussed above in (81)–(82).

**Step 5: Verification argument 1: the property of $\phi^*$**

In **Steps 5 and 6**, we verify that the guessed solution is the value function (18). **Steps 1 through 4** show that the control given in (38) is in equilibrium. The proof here follows the proof of Theorem 5.2 in Mbodji and Pirvu (2025)[55], with an adaptation to our setting because the studied system dynamics are different from each other.

Herein, we present several notations. The population driven by the control $\phi \in \mathbb{F}$ is denoted by $x$. Similarly, the populations driven by $\phi^* \in \mathbb{F}$ and $\phi_\varepsilon \in \mathbb{F}$ are denoted by $x^* \in B(0,+\infty)$ and $x^{(\varepsilon)} \in B(0,+\infty)$, respectively, which are bounded. For any $t \ge 0$, we have

$$\Delta_t\Phi\left(x_t^*\right)-\Phi\left(x_0\right)=\int_0^t\Delta_s\int_0^{+\infty}\left\{-r\phi^*(r)x_s^*(r)+R\left(1-x_s^*(r)\right)\right\}\Phi_f\left(r,x_s^*\right)p(\mathrm{d}r)\mathrm{d}s+\int_0^t\frac{\mathrm{d}\Delta_s}{\mathrm{d}s}\Phi\left(x_s^*\right)\mathrm{d}s\,. \tag{84}$$

Here, we used

$$\begin{aligned}\frac{\mathrm{d}}{\mathrm{d}t}\left(\Delta_t\Phi\left(x_t^*\right)\right)&=\frac{\mathrm{d}\Delta_t}{\mathrm{d}t}\Phi\left(x_t^*\right)+\Delta_t\frac{\mathrm{d}}{\mathrm{d}t}\Phi\left(x_t^*\right)\\&=\frac{\mathrm{d}\Delta_t}{\mathrm{d}t}\Phi\left(x_t^*\right)+\Delta_t\int_0^{+\infty}\frac{\partial x_t^*(r)}{\partial t}\Phi_f\left(r,x_s^*\right)p(\mathrm{d}r)\\&=\frac{\mathrm{d}\Delta_t}{\mathrm{d}t}\Phi\left(x_t^*\right)+\Delta_t\int_0^{+\infty}\left\{-r\phi^*(r)x_s^*(r)+R\left(1-x_s^*(r)\right)\right\}\Phi_f\left(r,x_s^*\right)p(\mathrm{d}r)\end{aligned}\tag{85}$$

and integrated it in time. From the extended HJB system, for any $f\in E$,

$$\begin{aligned}&\int_0^{+\infty}\left\{-r\phi^*(r)f(r)+R\left(1-f(r)\right)\right\}\Phi_f\left(r,f\right)p(\mathrm{d}r)\\&+\int_0^{+\infty}f(r)p(\mathrm{d}r)-\int_0^{+\infty}\frac{r}{\eta(r)}f(r)\left(\phi^*(r)\ln\phi^*(r)-\phi^*(r)+1\right)p(\mathrm{d}r)+\Psi_0(f)=0\end{aligned}\tag{86}$$

and

$$\begin{aligned}&\frac{\partial\Psi_t(f)}{\partial t}+\int_0^{+\infty}\left\{-r\phi^*(r)f(r)+R\left(1-f(r)\right)\right\}\Psi_{t,f}\left(r,f\right)p(\mathrm{d}r)\\&+\rho_t\left\{\int_0^{+\infty}f(r)p(\mathrm{d}r)-\int_0^{+\infty}\frac{r}{\eta(r)}f(r)\left(\phi^*(r)\ln\phi^*(r)-\phi^*(r)+1\right)p(\mathrm{d}r)\right\}\end{aligned}.\tag{87}$$

Then, we have

$$\begin{aligned}&\Delta_t\Phi\left(x_t^*\right)-\Phi\left(x_0\right)\\&=\int_0^t\Delta_s\int_0^{+\infty}\left\{-r\phi^*(r)x_s^*(r)+R\left(1-x_s^*(r)\right)\right\}\Phi_f\left(r,x_s^*\right)p(\mathrm{d}r)\mathrm{d}s+\int_0^t\frac{\mathrm{d}\Delta_s}{\mathrm{d}s}\Phi\left(x_s^*\right)\mathrm{d}s\\&=\int_0^t\Delta_s\left\{-\int_0^{+\infty}x_s^*(r)p(\mathrm{d}r)+\int_0^{+\infty}\frac{r}{\eta(r)}x_s^*(r)\left(\phi^*(r)\ln\phi^*(r)-\phi^*(r)+1\right)p(\mathrm{d}r)-\Psi_0\left(x_s^*\right)\right\}\mathrm{d}s\\&+\int_0^t\frac{\mathrm{d}\Delta_s}{\mathrm{d}s}\Phi\left(x_s^*\right)\mathrm{d}s\\&=\int_0^t\Delta_s\left\{-\int_0^{+\infty}x_s^*(r)p(\mathrm{d}r)+\int_0^{+\infty}\frac{r}{\eta(r)}x_s^*(r)\left(\phi^*(r)\ln\phi^*(r)-\phi^*(r)+1\right)p(\mathrm{d}r)\right\}\mathrm{d}s\\&-\int_0^t\left(\Delta_s\Psi_0\left(x_s^*\right)-\frac{\mathrm{d}\Delta_s}{\mathrm{d}s}\Phi\left(x_s^*\right)\right)\mathrm{d}s,\end{aligned}\tag{88}$$

where we used $f=x_s^*$ in (86) and (88). Subsequently, we obtain the following equation:

$$\begin{aligned}\Phi\left(x_0\right)=&\int_0^t\Delta_s\left\{\int_0^{+\infty}x_s^*(r)p(\mathrm{d}r)-\int_0^{+\infty}\frac{r}{\eta(r)}x_s^*(r)\left(\phi^*(r)\ln\phi^*(r)-\phi^*(r)+1\right)p(\mathrm{d}r)\right\}\mathrm{d}s\\&+\Delta_t\Phi\left(x_t^*\right)+\int_0^t\left(\Delta_s\Psi_0\left(x_s^*\right)-\frac{\mathrm{d}\Delta_s}{\mathrm{d}s}\Phi\left(x_s^*\right)\right)\mathrm{d}s\end{aligned}.\tag{89}$$

We show that the two terms in the second lines in (89) vanish at the limit $t\to+\infty$. First,

$$\lim_{t\to+\infty}\Delta_t=\lim_{t\to+\infty}\int_0^{+\infty}e^{-\delta t}\mu(\mathrm{d}\delta)=0,\tag{90}$$

showing that the first term in the second line in (89) vanishes when $t\to+\infty$. Second, we prove

$$\int_0^{+\infty}\left(\Delta_s\Psi_0\left(x_s^*\right)-\frac{\mathrm{d}\Delta_s}{\mathrm{d}s}\Phi\left(x_s^*\right)\right)\mathrm{d}s=0.\tag{91}$$

From (87) and (21), $\Psi_0$ admits the representation

$$\begin{aligned}
&\Psi_0\left(x_s^*\right)\\
&=\int_0^{+\infty}\int_0^{+\infty}\left(-\delta e^{-\delta\tau}\right)\left\{\int_0^{+\infty}x_{\tau+s}^*\left(r\right)p\left(\mathrm{d}r\right)-\int_0^{+\infty}\frac{r}{\eta\left(r\right)}x_{\tau+s}^*\left(r\right)\left(\phi^*\left(r\right)\ln\phi^*\left(r\right)-\phi^*\left(r\right)+1\right)p\left(\mathrm{d}r\right)\right\}\mathrm{d}\tau\mu\left(\mathrm{d}\delta\right) \quad (92)\\
&=\int_0^{+\infty}\frac{\mathrm{d}\Delta_\tau}{\mathrm{d}\tau}\left\{\int_0^{+\infty}x_{\tau+s}^*\left(r\right)p\left(\mathrm{d}r\right)-\int_0^{+\infty}\frac{r}{\eta\left(r\right)}x_{\tau+s}^*\left(r\right)\left(\phi^*\left(r\right)\ln\phi^*\left(r\right)-\phi^*\left(r\right)+1\right)p\left(\mathrm{d}r\right)\right\}\mathrm{d}\tau.
\end{aligned}$$

We also have

$$\Phi\left(x_s^*\right)=\int_0^{+\infty}\Delta_\tau\left\{\int_0^{+\infty}x_{\tau+s}^*\left(r\right)p\left(\mathrm{d}r\right)-\int_0^{+\infty}\frac{r}{\eta\left(r\right)}x_{\tau+s}^*\left(r\right)\left(\phi^*\left(r\right)\ln\phi^*\left(r\right)-\phi^*\left(r\right)+1\right)p\left(\mathrm{d}r\right)\right\}\mathrm{d}\tau. \quad (93)$$

We therefore obtain

$$\begin{aligned}
&\Delta_s\Psi_0\left(x_s^*\right)-\frac{\mathrm{d}\Delta_s}{\mathrm{d}s}\Phi\left(x_s^*\right)\\
&=\int_0^{+\infty}\Delta_s\frac{\mathrm{d}\Delta_\tau}{\mathrm{d}\tau}\left\{\int_0^{+\infty}x_{\tau+s}^*\left(r\right)p\left(\mathrm{d}r\right)-\int_0^{+\infty}\frac{r}{\eta\left(r\right)}x_{\tau+s}^*\left(r\right)\left(\phi^*\left(r\right)\ln\phi^*\left(r\right)-\phi^*\left(r\right)+1\right)p\left(\mathrm{d}r\right)\right\}\mathrm{d}\tau\\
&-\int_0^{+\infty}\frac{\mathrm{d}\Delta_s}{\mathrm{d}s}\Delta_\tau\left\{\int_0^{+\infty}x_{\tau+s}^*\left(r\right)p\left(\mathrm{d}r\right)-\int_0^{+\infty}\frac{r}{\eta\left(r\right)}x_{\tau+s}^*\left(r\right)\left(\phi^*\left(r\right)\ln\phi^*\left(r\right)-\phi^*\left(r\right)+1\right)p\left(\mathrm{d}r\right)\right\}\mathrm{d}\tau \quad (94)\\
&=\int_0^{+\infty}\left(\Delta_s\frac{\mathrm{d}\Delta_\tau}{\mathrm{d}\tau}-\frac{\mathrm{d}\Delta_s}{\mathrm{d}s}\Delta_\tau\right)\left\{\begin{aligned}&\int_0^{+\infty}x_{\tau+s}^*\left(r\right)p\left(\mathrm{d}r\right)\\&-\int_0^{+\infty}\frac{r}{\eta\left(r\right)}x_{\tau+s}^*\left(r\right)\left(\phi^*\left(r\right)\ln\phi^*\left(r\right)-\phi^*\left(r\right)+1\right)p\left(\mathrm{d}r\right)\end{aligned}\right\}\mathrm{d}\tau\\
&=0
\end{aligned}$$

due to the symmetry in $\Delta_s\frac{\mathrm{d}\Delta_\tau}{\mathrm{d}\tau}-\frac{\mathrm{d}\Delta_s}{\mathrm{d}s}\Delta_\tau$, proving (91). Consequently, it follows that

$$\begin{aligned}
\Phi\left(x_0\right)&=\int_0^{+\infty}\Delta_s\left\{\int_0^{+\infty}x_s^*\left(r\right)p\left(\mathrm{d}r\right)-\int_0^{+\infty}\frac{r}{\eta\left(r\right)}x_s^*\left(r\right)\left(\phi^*\left(r\right)\ln\phi^*\left(r\right)-\phi^*\left(r\right)+1\right)p\left(\mathrm{d}r\right)\right\}\mathrm{d}s\\
&=\theta\left(x_0;\phi^*\right)
\end{aligned}. \quad (95)$$

**Step 6: Verification argument 2: Equilibrium condition**

The last step is to show (38) and (39). For later use, for any $f\in E$ and $\varphi\in\mathbb{F}$, we set

$$\mathbb{H}\left(f,\varphi\right)=\int_0^{+\infty}f\left(r\right)p\left(\mathrm{d}r\right)-\int_0^{+\infty}\frac{r}{\eta\left(r\right)}f\left(r\right)\left(\varphi\left(r\right)\ln\varphi\left(r\right)-\varphi\left(r\right)+1\right)p\left(\mathrm{d}r\right). \quad (96)$$

Because $\phi_\varepsilon=\phi^*$ for $\varepsilon\le s$ and $\phi_\varepsilon=\phi$ for $0\le s<\varepsilon$, we have

$$
\begin{aligned}
&\frac{\theta\left(x;\phi^*\right)-\theta\left(x;\phi_\varepsilon\right)}{\varepsilon}\\
&=\frac{1}{\varepsilon}\int_0^{+\infty}\Delta_s\left(\mathbb{H}\left(x_s^*,\phi^*\right)-\mathbb{H}\left(x_s^{(\varepsilon)},\phi_{\varepsilon,s}\right)\right)\mathrm{d}s\\
&=\frac{1}{\varepsilon}\int_0^{\varepsilon}\Delta_s\left(\mathbb{H}\left(x_s^*,\phi^*\right)-\mathbb{H}\left(x_s^{(\varepsilon)},\phi_s\right)\right)\mathrm{d}s\\
&+\frac{1}{\varepsilon}\int_\varepsilon^{+\infty}\left(\Delta_{s-\varepsilon}-\Delta_s\right)\left(\mathbb{H}\left(x_s^{(\varepsilon)},\phi^*\right)-\mathbb{H}\left(x_s^*,\phi^*\right)\right)\mathrm{d}s+\frac{1}{\varepsilon}\int_\varepsilon^{+\infty}\Delta_{s-\varepsilon}\left(\mathbb{H}\left(x_s^*,\phi^*\right)-\mathbb{H}\left(x_s^{(\varepsilon)},\phi^*\right)\right)\mathrm{d}s\\
&=K_1\left(\varepsilon\right)+K_2\left(\varepsilon\right)+K_3\left(\varepsilon\right),
\end{aligned}
\tag{97}
$$

where

$$
K_1\left(\varepsilon\right)=\frac{1}{\varepsilon}\int_0^{\varepsilon}\Delta_s\left(\mathbb{H}\left(x_s^*,\phi^*\right)-\mathbb{H}\left(x_s^{(\varepsilon)},\phi_s\right)\right)\mathrm{d}s, \tag{98}
$$

$$
K_2\left(\varepsilon\right)=\frac{1}{\varepsilon}\int_\varepsilon^{+\infty}\left(\Delta_{s-\varepsilon}-\Delta_s\right)\left(\mathbb{H}\left(x_s^{(\varepsilon)},\phi^*\right)-\mathbb{H}\left(x_s^*,\phi^*\right)\right)\mathrm{d}s, \tag{99}
$$

and

$$
K_3\left(\varepsilon\right)=\frac{1}{\varepsilon}\int_\varepsilon^{+\infty}\Delta_{s-\varepsilon}\left(\mathbb{H}\left(x_s^*,\phi^*\right)-\mathbb{H}\left(x_s^{(\varepsilon)},\phi^*\right)\right)\mathrm{d}s. \tag{100}
$$

We evaluate each of (98)–(100).

First, for $K_1\left(\varepsilon\right)$, by the mean value theorem, with some $\phi_0>0$, we have

$$
\lim_{\varepsilon\to 0}K_1\left(\varepsilon\right)=\lim_{\varepsilon\to 0}\frac{1}{\varepsilon}\int_0^{\varepsilon}\Delta_s\left(\mathbb{H}\left(x_s^*,\phi^*\right)-\mathbb{H}\left(x_s^{(\varepsilon)},\phi_s\right)\right)\mathrm{d}s=\mathbb{H}\left(x_0,\phi^*\right)-\mathbb{H}\left(x_0,\phi_0\right). \tag{101}
$$

Second, for $K_2\left(\varepsilon\right)$, we have (we use (56))

$$
\begin{aligned}
\left|K_2\left(\varepsilon\right)\right|&\le-\int_\varepsilon^{+\infty}\frac{\Delta_s-\Delta_{s-\varepsilon}}{\varepsilon}\left|\mathbb{H}\left(x_s^{(\varepsilon)},\phi^*\right)-\mathbb{H}\left(x_s^*,\phi^*\right)\right|\mathrm{d}s\\
&=-\int_\varepsilon^{+\infty}\frac{\mathrm{d}\Delta_{s_\varepsilon}}{\mathrm{d}s}\left|\mathbb{H}\left(x_s^{(\varepsilon)},\phi^*\right)-\mathbb{H}\left(x_s^*,\phi^*\right)\right|\mathrm{d}s\\
&=-\int_\varepsilon^{+\infty}\frac{\mathrm{d}\Delta_{s_\varepsilon}}{\mathrm{d}s}\left|\begin{gathered}\int_0^{+\infty}\left(x_s^{(\varepsilon)}\left(r\right)-x_s^*\left(r\right)\right)p\left(\mathrm{d}r\right)\\-\int_0^{+\infty}\frac{r}{\eta\left(r\right)}\left(x_s^{(\varepsilon)}\left(r\right)-x_s^*\left(r\right)\right)\left(\phi^*\left(r\right)\ln\phi^*\left(r\right)^*-\phi^*\left(r\right)+1\right)p\left(\mathrm{d}r\right)\end{gathered}\right|\mathrm{d}s\\
&=-\int_\varepsilon^{+\infty}\frac{\mathrm{d}\Delta_{s_\varepsilon}}{\mathrm{d}s}\cdot 2\left|\int_0^{+\infty}p\left(\mathrm{d}r\right)+\int_0^{+\infty}\frac{r}{\eta\left(r\right)}\left(\phi^*\left(r\right)\ln\phi^*\left(r\right)^*-\phi^*\left(r\right)+1\right)p\left(\mathrm{d}r\right)\right|\mathrm{d}s\\
&\le 4\int_\varepsilon^{+\infty}\left(-\frac{\mathrm{d}\Delta_{s_\varepsilon}}{\mathrm{d}s}\right)\mathrm{d}s
\end{aligned}
\tag{102}
$$

with some $s_\varepsilon\in\left(s-\varepsilon,s\right)$ by the continuity of $\Delta_s$ ($s>0$), and we used (96) and that $\phi^*$ is bounded between 0 and 1. Because $\Delta_s$ is not increasing and convex, the last integral is evaluated as

$$
0\le\int_\varepsilon^{+\infty}\left(-\frac{\mathrm{d}\Delta_{s_\varepsilon}}{\mathrm{d}s}\right)\mathrm{d}s\le\int_\varepsilon^{+\infty}\left(-\frac{\mathrm{d}\Delta_{s-\varepsilon}}{\mathrm{d}s}\right)\mathrm{d}s=\int_0^{+\infty}\left(-\frac{\mathrm{d}\Delta_s}{\mathrm{d}s}\right)\mathrm{d}s=\Delta_0=1, \tag{103}
$$

showing that $K_2(\varepsilon)$ is uniformly bounded irrespective of $\varepsilon > 0$. Then, by dominated convergence in conjunction with $\lim_{\varepsilon \to 0} \left| \mathbb{H}\left(x_s^{(\varepsilon)}, \phi^*\right) - \mathbb{H}\left(x_s^*, \phi^*\right) \right| = 0$ at each $s > 0$ (this follows from the fact that $\phi_{\varepsilon,s} \to \phi^*$ as $\varepsilon \to 0$ at each $s > 0$ and the form of the solution (10)), we obtain

$$
\begin{aligned}
\left|K_2(\varepsilon)\right| &= -\int_{\varepsilon}^{+\infty} \frac{\mathrm{d}\Delta_{s_\varepsilon}}{\mathrm{d}s} \left| \mathbb{H}\left(x_s^{(\varepsilon)}, \phi^*\right) - \mathbb{H}\left(x_s^*, \phi^*\right) \right| \mathrm{d}s \\
&= -\int_{\varepsilon}^{+\infty} \frac{\mathrm{d}\Delta_{s-\varepsilon}}{\mathrm{d}s} \left| \mathbb{H}\left(x_s^{(\varepsilon)}, \phi^*\right) - \mathbb{H}\left(x_s^*, \phi^*\right) \right| \mathrm{d}s \\
&\to -\int_{0}^{+\infty} \frac{\mathrm{d}\Delta_s}{\mathrm{d}s} \cdot 0 \mathrm{d}s \\
&= 0.
\end{aligned}
\tag{104}
$$

Third, for $K_3(\varepsilon)$, we have

$$
\begin{aligned}
K_3(\varepsilon) &= \frac{1}{\varepsilon} \int_{\varepsilon}^{+\infty} \Delta_{s-\varepsilon} \left( \mathbb{H}\left(x_s^*, \phi^*\right) - \mathbb{H}\left(x_s^{(\varepsilon)}, \phi^*\right) \right) \mathrm{d}s \\
&= \frac{1}{\varepsilon} \int_{0}^{+\infty} \Delta_s \left( \mathbb{H}\left(x_{s+\varepsilon}^*, \phi^*\right) - \mathbb{H}\left(x_{s+\varepsilon}^{(\varepsilon)}, \phi^*\right) \right) \mathrm{d}s \\
&= \frac{1}{\varepsilon} \int_{0}^{+\infty} \Delta_s \left( \mathbb{H}\left(x_{s+\varepsilon}^*, \phi^*\right) - \mathbb{H}\left(x_s^*, \phi^*\right) \right) \mathrm{d}s + \frac{1}{\varepsilon} \int_{0}^{+\infty} \Delta_s \left( \mathbb{H}\left(x_s^*, \phi^*\right) - \mathbb{H}\left(x_{s+\varepsilon}^{(\varepsilon)}, \phi^*\right) \right) \mathrm{d}s \\
&= \frac{1}{\varepsilon} \left( \theta\left(x_\varepsilon^*; \phi^*\right) - \theta\left(x_0; \phi^*\right) \right) + \frac{1}{\varepsilon} \left( \theta\left(x_0; \phi^*\right) - \theta\left(x_\varepsilon^{(\varepsilon)}; \phi^*\right) \right) \\
&= \frac{1}{\varepsilon} \left( \Phi\left(x_\varepsilon^*\right) - \Phi\left(x_0\right) \right) + \frac{1}{\varepsilon} \left( \Phi\left(x_0\right) - \Phi\left(x_\varepsilon^{(\varepsilon)}\right) \right)
\end{aligned}
\tag{105}
$$

and

$$
\begin{aligned}
\frac{1}{\varepsilon} \left( \Phi\left(x_\varepsilon^*\right) - \Phi\left(x_0\right) \right) &= \frac{1}{\varepsilon} \int_0^{\varepsilon} \int_0^{+\infty} \left\{ -r\phi^*(r) x_s^*(r) + R\left(1 - x_s^*(r)\right) \right\} \Phi_f\left(r, x_s^*\right) p(\mathrm{d}r) \mathrm{d}s \\
&\underset{\varepsilon \to 0}{\to} \int_0^{+\infty} \left\{ -r\phi^*(r) x_0(r) + R\left(1 - x_0(r)\right) \right\} \Phi_f\left(r, x_0\right) p(\mathrm{d}r)
\end{aligned}
\tag{106}
$$

as well as

$$
\begin{aligned}
\frac{1}{\varepsilon} \left( \Phi\left(x_0\right) - \Phi\left(x_\varepsilon^{(\varepsilon)}\right) \right) &= -\frac{1}{\varepsilon} \int_0^{\varepsilon} \int_0^{+\infty} \left\{ -r\phi^*(r) x_{\varepsilon+s}^{(\varepsilon)}(r) + R\left(1 - x_{\varepsilon+s}^{(\varepsilon)}(r)\right) \right\} \Phi_f\left(r, x_s^*\right) p(\mathrm{d}r) \mathrm{d}s \\
&\underset{\varepsilon \to 0}{\to} - \int_0^{+\infty} \left\{ -r\phi(r) x_0(r) + R\left(1 - x_0(r)\right) \right\} \Phi_f\left(r, x_0\right) p(\mathrm{d}r)
\end{aligned}
\tag{107}
$$

because $x_\varepsilon^{(\varepsilon)} = x_0$. Consequently, we arrive at the desired result

$$
\begin{aligned}
\liminf_{\varepsilon\to+0}\frac{\theta\left(x;\phi^*\right)-\theta\left(x;\phi_\varepsilon\right)}{\varepsilon} &= \liminf_{\varepsilon\to+0}\left\{K_1\left(\varepsilon\right)+K_2\left(\varepsilon\right)+K_3\left(\varepsilon\right)\right\} \\
&= \lim_{\varepsilon\to+0}\left\{K_1\left(\varepsilon\right)+K_2\left(\varepsilon\right)+K_3\left(\varepsilon\right)\right\} \\
&= \mathbb{H}\left(x_0,\phi^*\right)-\mathbb{H}\left(x_0,\phi_0\right) \\
&\quad +\int_0^{+\infty}\left\{-r\phi^*\left(r\right)x_0\left(r\right)+R\left(1-x_0\left(r\right)\right)\right\}\Phi_f\left(r,x_0\right)p\left(\mathrm{d}r\right) \\
&\quad -\int_0^{+\infty}\left\{-r\phi_0\left(r\right)x_0\left(r\right)+R\left(1-x_0\left(r\right)\right)\right\}\Phi_f\left(r,x_0\right)p\left(\mathrm{d}r\right) \\
&= \mathbb{H}\left(x_0,\phi^*\right)+\int_0^{+\infty}\left\{-r\phi^*\left(r\right)x_0\left(r\right)+R\left(1-x_0\left(r\right)\right)\right\}\Phi_f\left(r,x_0\right)p\left(\mathrm{d}r\right) \\
&\quad -\left\{\mathbb{H}\left(x_0,\phi_0\right)+\int_0^{+\infty}\left\{-r\phi_0\left(r\right)x_0\left(r\right)+R\left(1-x_0\left(r\right)\right)\right\}\Phi_f\left(r,x_0\right)p\left(\mathrm{d}r\right)\right\} \\
&= \sup_{\phi}\left\{\mathbb{H}\left(x_0,\phi\right)+\int_0^{+\infty}\left\{-r\phi\left(r\right)x_0\left(r\right)+R\left(1-x_0\left(r\right)\right)\right\}\Phi_f\left(r,x_0\right)p\left(\mathrm{d}r\right)\right\} \\
&\quad -\left\{\mathbb{H}\left(x_0,\phi_0\right)+\int_0^{+\infty}\left\{-r\phi_0\left(r\right)x_0\left(r\right)+R\left(1-x_0\left(r\right)\right)\right\}\Phi_f\left(r,x_0\right)p\left(\mathrm{d}r\right)\right\} \\
&\geq 0.
\end{aligned}
\tag{108}
$$

Moreover, this shows that $\phi^*$ is an equilibrium control.

□

### A2. Finite-dimensional system

### A2.1 Representation via quantization

In the main text, we presented a finite-dimensional version of the extended HJB system, which is a quantized version of the infinite-dimensional system studied in the main text. We assume that $N$ samples $(r_1, r_2, ..., r_N)$ and $M$ samples $(\delta_1, \delta_2, ..., \delta_M)$ were drawn from the probability measures $p$ and $\mu$, respectively. The samples can be either random, as in classical Monte Carlo methods, or deterministic, as in the quantile-based quantization method (e.g., Yoshioka and Yoshooka, 2024b)[33]. For binary variables $x_i \in \{0,1\}$ ($i = 1,2,3,...,N$), we write $x = (x_1, x_2, ..., x_N)$. Note that the convergence argument in this section is only formal, which is why we addressed the numerical convergence experiment in **Section A2.2**; nevertheless, the results in this section would be useful for better understanding the infinite-dimensional model discussed in the main text.

For numerical computation of the extended HJB system, we set each $r_i$ as

$$\int_0^{r_i} p(\mathrm{d}r) = \frac{2i-1}{2N}, \quad i = 1,2,3,...,N. \tag{109}$$

Similarly, we set each $\delta_j$ as

$$\int_0^{\delta_j} \mu(\mathrm{d}\delta) = \frac{2j-1}{2M}, \quad j = 1,2,3,...,M. \tag{110}$$

With this quantization method, we can empirically approximate the two probability measures $p$ and $\mu$ with maximum errors of $O(N^{-1})$ and $O(M^{-1})$, respectively, in terms of the empirical and true cumulative distribution functions (e.g., Proof of Proposition 2 in Yoshioka [78]).

Under these preparations, the extended HJB system in the finite-dimensional setting reads

$$\sup_{\phi} \left\{ \begin{array}{l} -\sum_{i=1}^{N} r_i x_i \phi_i \left(\Phi - \bar{\Phi}_i\right) - \sum_{i=1}^{N} R\left(1 - x_i\right)\left(\Phi - \hat{\Phi}_i\right) + \frac{1}{N}\sum_{i=1}^{N} x_i \\ -\frac{1}{N}\sum_{i=1}^{N} \frac{1}{\eta} r_i x_i \left(\phi_i \ln \phi_i - \phi_i + 1\right) + \Psi(0,x) \end{array} \right\} = 0 \tag{111}$$

with

$$\bar{\Phi}_i = \Phi\Big(x_1, x_2, ..., \underset{i}{0}, ...\Big) \text{ and } \hat{\Phi}_i = \Phi\Big(x_1, x_2, ..., \underset{i}{1}, ...\Big) \tag{112}$$

and

$$\phi_i^* = e^{-\eta\left(\Phi - \bar{\Phi}_i\right)}. \tag{113}$$

Here, $\Psi(t,x)$ satisfies

$$\begin{array}{l} \frac{\partial \Psi(t,x)}{\partial t} - \sum_{i=1}^{N} r_i \phi_i^* x_i \left(\Psi - \bar{\Psi}_i\right) - \sum_{i=1}^{N} R\left(1 - x_i\right)\left(\Psi - \hat{\Psi}_i\right) \\ + \left(\frac{\mathrm{d}}{\mathrm{d}t}\Delta_t\right)\left\{\frac{1}{N}\sum_{i=1}^{N} x_i - \frac{1}{N}\sum_{i=1}^{N}\frac{1}{\eta} r_i x_i \left(\phi_i^* \ln \phi_i^* - \phi_i^* + 1\right)\right\} = 0 \end{array}, \quad t > 0 \tag{114}$$

with $\Psi(+\infty, x) = 0$, where

$$\Delta_t = \frac{1}{M}\sum_{j=1}^{M} e^{-\delta_j t} \quad \text{and} \quad \frac{\mathrm{d}}{\mathrm{d}t}\Delta_t = -\frac{1}{M}\sum_{j=1}^{M} c_j e^{-\delta_j t} . \tag{115}$$

For **Proposition** 1, we guess the solution of the form

$$\Phi = \frac{1}{M}\sum_{i=1}^{N} A_i x_i + B \quad \text{and} \quad \Psi = \frac{1}{M}\sum_{j=1}^{M}\left(\frac{1}{N}\sum_{i=1}^{N} C_{i,j} x_i + D_j\right) e^{-\delta_j t} , \tag{116}$$

with which we obtain

$$\phi_i^* = e^{-\eta A_i} . \tag{117}$$

The equations (111) and (114) reduce to

$$-\frac{1}{N}\sum_{i=1}^{N} r_i x_i \frac{1-e^{-\eta A_i}}{\eta} + \frac{1}{N}\sum_{i=1}^{N} R(1-x_i)A_i + \frac{1}{N}\sum_{i=1}^{N} x_i + \frac{1}{N}\sum_{i=1}^{N}\left(\frac{1}{M}\sum_{j=1}^{M} C_{i,j}\right)x_i + \frac{1}{M}\sum_{j=1}^{M} D_j = 0 \tag{118}$$

and

$$\begin{aligned}
&-\frac{1}{M}\sum_{j=1}^{M}\delta_j\left(\frac{1}{N}\sum_{i=1}^{N} C_{i,j} x_i + D_j\right)e^{-\delta_j t} \\
&-\frac{1}{N}\sum_{i=1}^{N} r_i \phi_i^* x_i\left(\frac{1}{M}\sum_{j=1}^{M} C_{i,j}\right)e^{-\delta_j t} + \frac{1}{N}\sum_{i=1}^{N} R(1-x_i)\left(\frac{1}{M}\sum_{j=1}^{M} C_{i,j}\right)e^{-\delta_j t} \\
&-\frac{1}{M}\sum_{j=1}^{M}\delta_j e^{-\delta_j t}\frac{1}{N}\sum_{i=1}^{N} x_i\left(1-\frac{1}{\eta}r_i\left(\phi_i^*\ln\phi_i^* - \phi_i^* + 1\right)\right) = 0,
\end{aligned} \tag{119}$$

respectively. From (118) and (119), we obtain the following system of algebraic equations:

$$A_i = \frac{1}{R}\left(-r_i\frac{1-e^{-\eta A_i}}{\eta} + + \frac{1}{M}\sum_{j=1}^{M} C_{i,j}\right), \tag{120}$$

$$R\frac{1}{N}\sum_{i=1}^{N} A_i + \frac{1}{M}\sum_{j=1}^{M} D_j = 0 , \tag{121}$$

$$C_{i,j} = -\frac{1}{MN}\frac{\delta_j}{\delta_j + r_i e^{-\eta A_i} + R}\left(1-\frac{r_i}{\eta}\left(1-\left(1+\eta A_i\right)e^{-\eta A_i}\right)\right), \tag{122}$$

and

$$-\delta_j D_j + R\frac{1}{N}\sum_{i=1}^{N} C_{i,j} = 0 . \tag{123}$$

The equation to find $A_i$ is given by

$$A_i = \frac{1}{R}\left(-\frac{1}{N}r_i\frac{1-e^{-\eta A_i}}{\eta} + \frac{1}{N} - \frac{1}{N}\left(1-\frac{r_i}{\eta}\left(1-\left(1+\eta A_i\right)e^{-\eta A_i}\right)\right)\frac{1}{M}\sum_{j=1}^{M}\frac{\delta_j}{\delta_j + r_i e^{-\eta A_i} + R}\right), \tag{124}$$

with which $C_{i,j}$ is obtained from (122), and $D_j$ from (123).

We can formally take the infinite-dimensional limit by taking $N, M \to +\infty$ :

$$A_i \to A(r) \quad \text{and} \quad C_{i,j} \to C(r,\delta) , \tag{125}$$

with which we obtain the corresponding infinite-dimensional results (33)–(36) in **Proposition 1**. From a computational perspective, the finite-dimensional system (120)–(124) would serve as an approximation of the infinite-dimensional system. The optimality result in **Proposition 1** applies to the finite-dimensional case because this case corresponds to the use of a discrete $p, \mu$.

A theoretical difference between finite-dimensional and infinite-dimensional systems is that the regularity conditions of the probability measures, $p, \mu$ (i.e., (31)-(32)), are necessary only for the latter. These conditions arise from the singularity of $p, \mu$ near the origin, which does not appear in a finite-dimensional setting. Another theoretical difference between them is in the extended HJB systems, where the finite-dimensional system is nonlocal, whereas the infinite-dimensional system is local. This difference arises from the corresponding system dynamics; the dynamics are driven by random jumps in the finite-dimensional system, while they are deterministic, and no jump is present in the infinite-dimensional system. We conjecture that this correspondence carries over to more complex population dynamics, although this is beyond the scope of this study.

### A2.2 Convergence test

We conclude this appendix by reporting on a convergence study of a finite-dimensional system. We regard the computed coefficient $B$ of the value function $V$, which is $\Phi$ according to **Proposition 1**, with the resolution $N = M = 4096$ as the reference solution because the integrals constituting $B$ do not seem to have exact formulae. We track the absolute errors between the reference and coarser numerical solutions and demonstrate that the choice of $N = M = 1024$ in the main text is reasonable.

**Table A1** lists the maximum error defined by $\max_{k=1,2,3,...,100} \left| B_{\mathrm{ref},k} - B_{\mathrm{num},k} \right|$, where the subscripts "ref" and "num" represent reference and coarser numerical solutions, respectively, and the other subscript $k$ represents the uncertainty-aversion coefficient $c = (k - 1/2)/100$ ($k = 1,2,3,...,100$). Therefore, the error quantifies the accuracy of the computed $B$ for a broad range of $c$ values. We select $\alpha_\delta = 1.25$ and $R = 0.02$, and the results are similar for the other parameter values examined in the main text. The reference $B$ values range from 53.4692 ($k = 1$) to 63.6110 ($k = 100$). **Table A1** suggests that the maximum error is at least halved as the resolution is doubled and is sufficiently small at $N = M = 1024$, where the maximum error relative to the reference is less than 0.01%. This result, combined with the error in the finer-resolution $N = M = 2048$, implies that the maximum error between the numerical solution at $N = M = 1024$ and the true $B$ value is at most $O(10^{-2})$, which is sufficiently small for the computational study described in **Section 4**.

**Table A1.** Maximum errors between reference and numerical $B$ .

| Resolution | Maximum error |
|---|---|
| $N = M = 256$ | 0.0127 |
| $N = M = 512$ | 0.0060 |
| $N = M = 1024$ | 0.0026 |
| $N = M = 2048$ | 0.0009 |

## A3. Experimental setting

This section briefly explains the experimental setting used to find the parameter values of $\alpha_r$ and $\beta_r$ in **Section 4**. The explanation here is adapted from Yoshioka and Hamagami (2024)[38]. The experimental flume presented in the middle of **Figure 1** was placed at the Faculty of Agriculture, Iwate University, Japan. The width and slope of the flume were 0.15 (m) and 1/70, respectively. We placed hemispheres with a diameter of 0.075 (m) covered by *Cladophora glomerata* in the observable area of the flume. After establishing a stationary flow with a discharge of 0.008 ($m^3$/s), we supplied sand particles with a unit-width sediment discharge of 0.00004 ($m^2$/s). The mean diameter of the sand particles supplied to the water flow was 2 mm based on one hundred samples with diameters between 1.70 mm and 2.36 mm. The algae population, which is the covering ratio of each gravel, was subsequently measured at 1, 2, 3, 4, 5, and 6 hours by a video camera (Victor: GZ-HD40, frame rate of 30 fps), from which we could evaluate $X$. Then, a least-squares fitting of the observed and experimental $X$ was conducted by ignoring the influence of growth to find values of $\alpha_r$ and $\beta_r$ (see **Figure A1**).

## A4. Auxiliary computational results

We present the computational results with $R = 0.005$ (1/day) and $R = 0.04$ (1/day) to show that the results and discussion in **Section 4** are not critically affected by the chosen $R$ values. **Figure A2** compares the computed $\phi^*$ for $\alpha_\delta = 1.25$ with $R = 0.005$ and $R = 0.04$. **Figure A3** compares the computed $\text{Ren-Env}$ relationships for $\alpha_\delta = 1.25$ and different values of $R$. The plots for each $R$ are similar, and a larger value of $R$ leads to more pessimistic results with a larger environmental index and model uncertainties, as discussed in **Section 4**.

We also present additional computational results where different $\alpha_r$ values are examined. **Figure A4** shows $V(1)$ values for different values of $\alpha_r = 0.295m$ ($m = 1,\ 2,\ 4,\ 100$) with $\alpha_\delta = 2$ being fixed. Here, we fix $\alpha_r \beta_r$ so that the average of $r$ does not change, so the case $m = 100$ is closest to an exponential case with a decay rate of $\alpha_r \beta_r$. **Figure A4** shows that $V(1)$ decreases as $\alpha_r$ decreases, i.e., as $p$ approaches exponential distributions. The model with $p$ closer to be exponential thus predicts smaller $V(1)$, namely rapid decay of algae population. This seems to yield generally smaller Ren values as shown in **Figure A5**. Although not presented here, a (nearly) exponential model with more rapid decay gives smaller, more optimistic $V(1)$, suggesting that considering the long-memory decay is crucial for evaluating the environmental index.

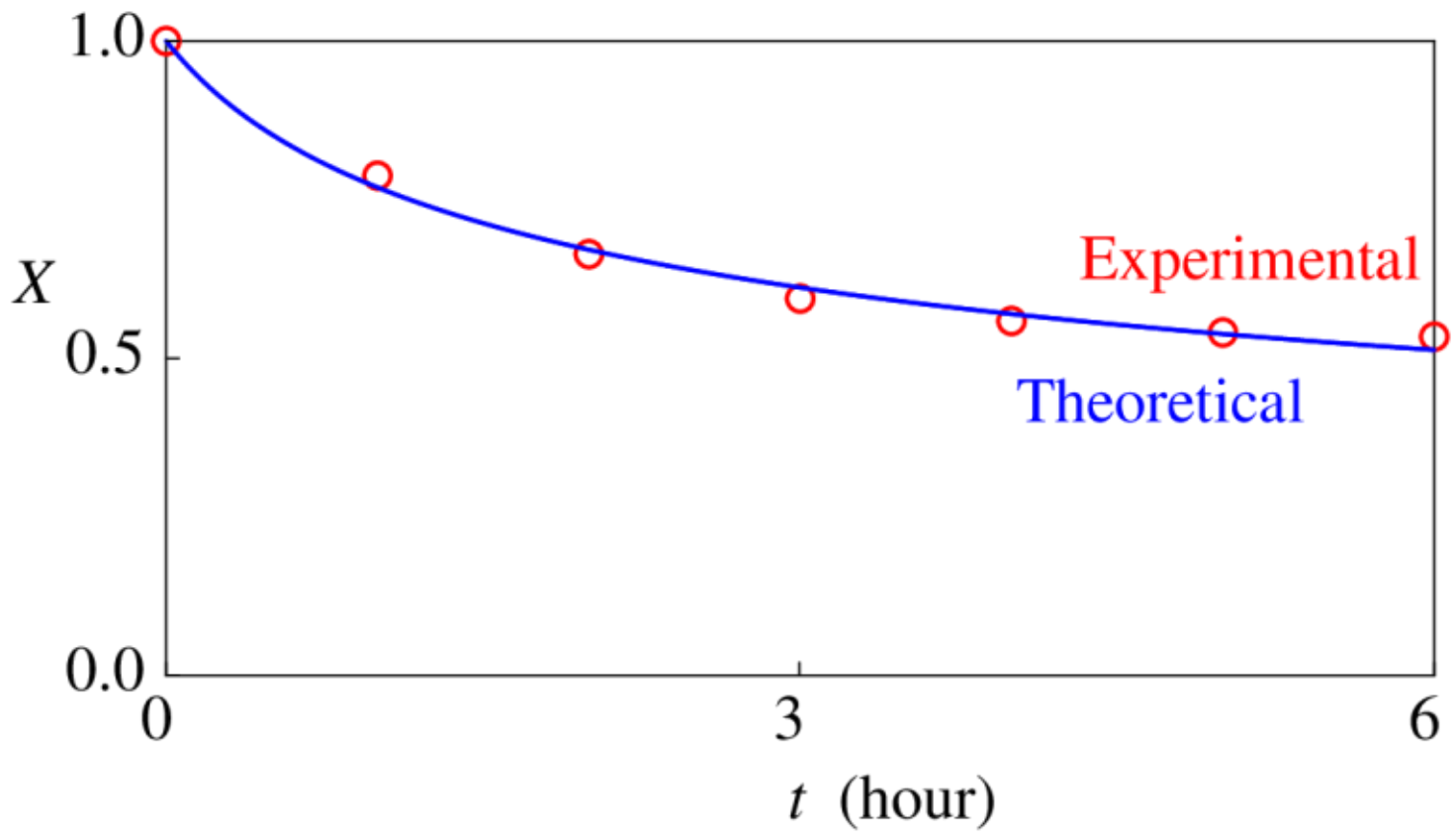

**Figure A1.** Empirical and theoretical algae cover for the model used in **Section 4**.

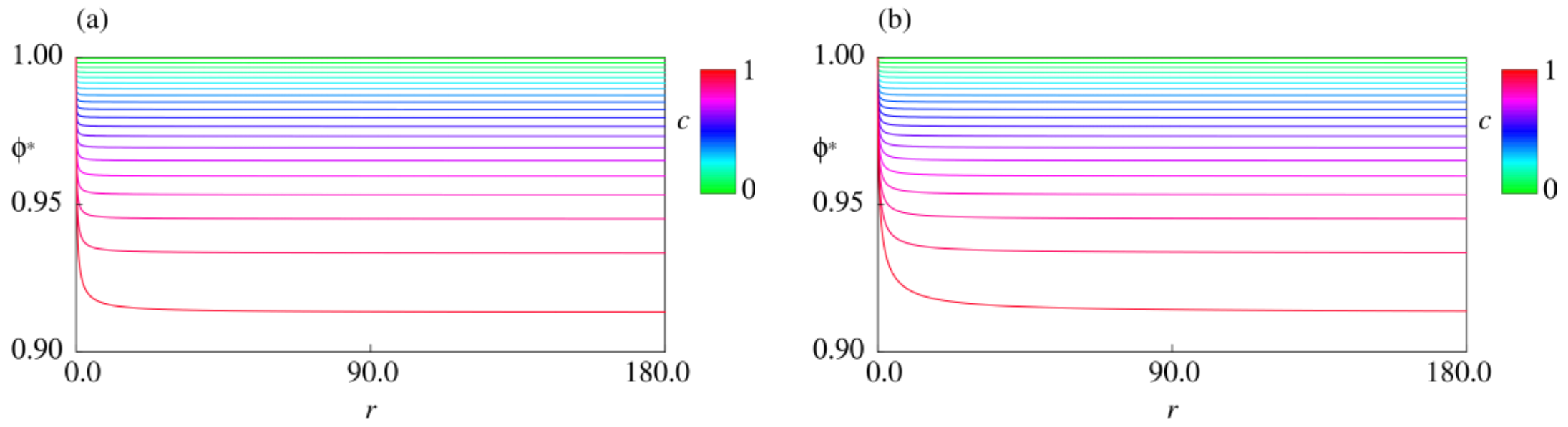

**Figure A2.** Computed $\phi^{*}$ for $\alpha_{\delta} = 1.25$ and different values of $R$ : (a) 0.005 and (b) 0.04.

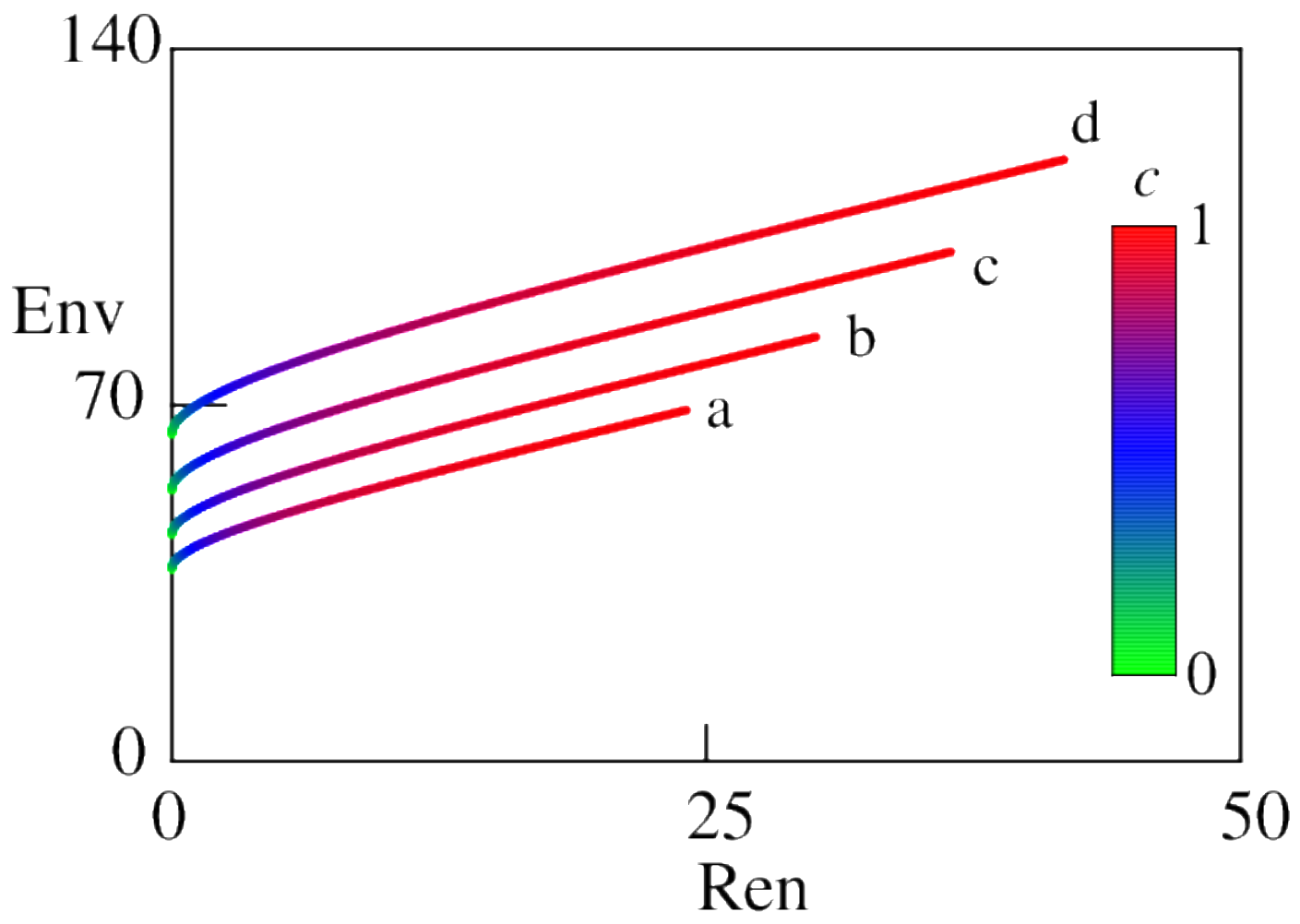

**Figure A3.** Computed Ren-Env relationship for $\alpha_{\delta} = 1.25$ and different values of $R$ : (a) 0.005, (b) 0.01, (c) 0.02, and (d) 0.04.

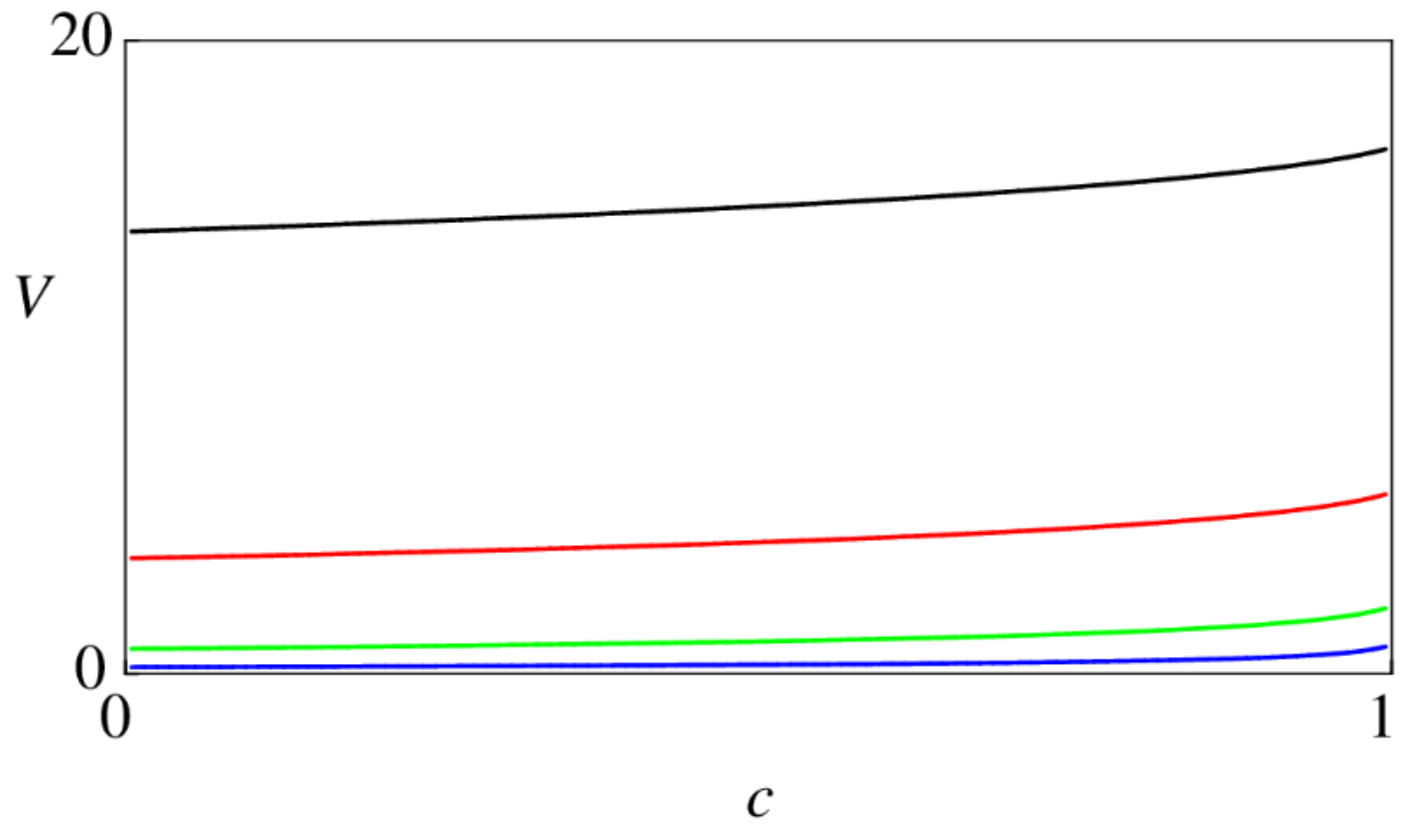

**Figure A4.** Computed $V = V(1)$ for $R = 0.01$ and different values of $\alpha_r = 0.295m$: $m = 1$ (black), $m = 2$ (red), $m = 4$ (green), and $m = 100$ (blue).

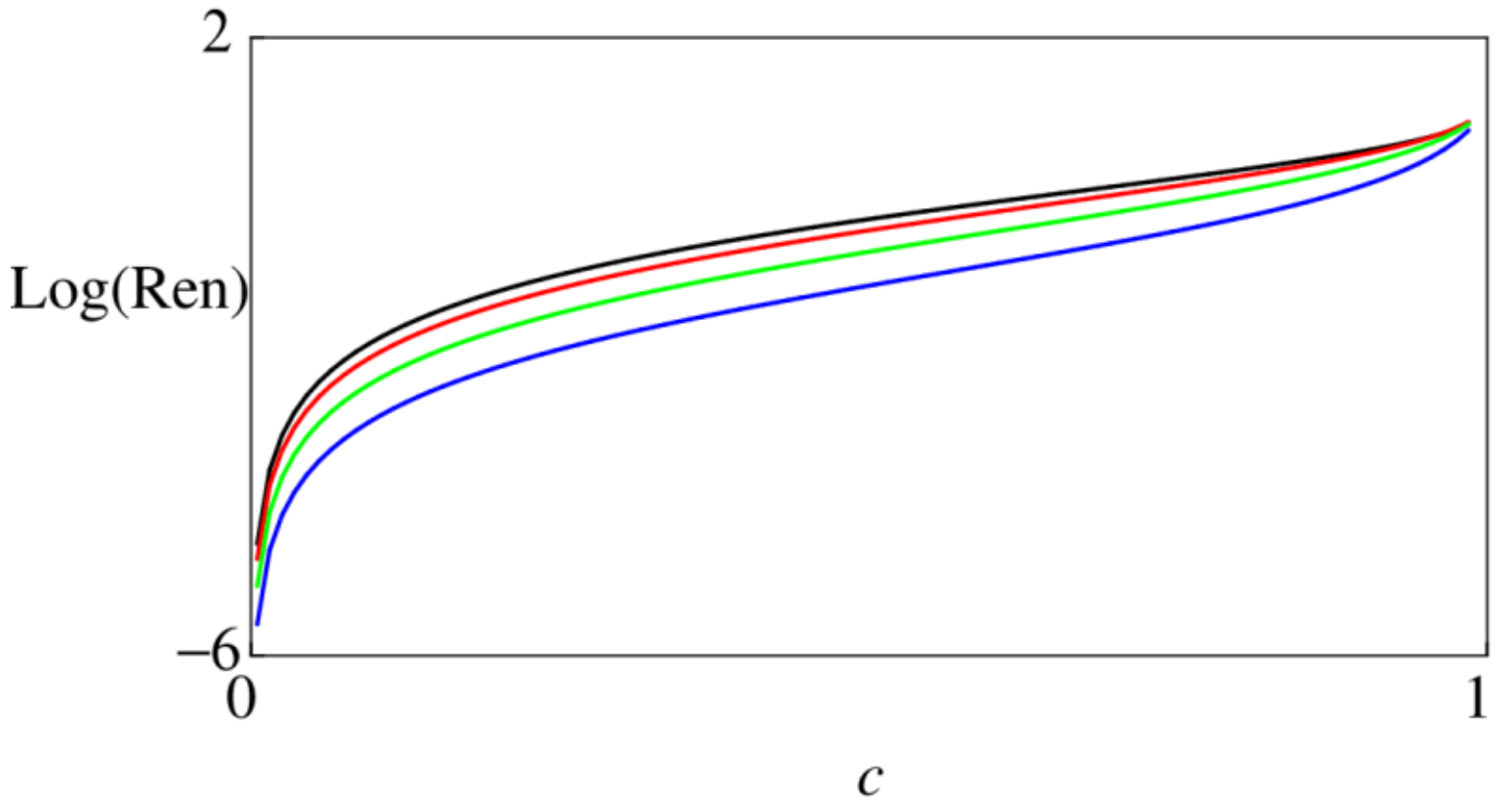

**Figure A5.** Computed (ordinary logarithm of) Ren in for $R = 0.01$ and different values of $\alpha_r = 0.295m$: $m = 1$ (black), $m = 2$ (red), $m = 4$ (green), and $m = 100$ (blue).

**References**

[1] Fritsch, M., Haupt, H., Flock, T., Schnurbus, J., & Sibbertsen, P. (2025). The memory puzzle in precipitation: Uncertainty in memory parameter estimation and implications for forecasting practice. International Journal of Forecasting. Online published. https://doi.org/10.1016/j.ijforecast.2025.11.003

[2] Granata, F., & Di Nunno, F. (2025). The anatomy of drought in Italy: statistical signatures, spatiotemporal persistence, and forecasting potential. Journal of Hydrology, 134428. https://doi.org/10.1016/j.jhydrol.2025.134428

[3] Shi, H., Yang, X., Tang, H., & Tu, Y. (2025). Dynamic prediction of PM2.5 concentration in China using experience replay with multiperiod memory buffers. Atmospheric Research, 320, 108063. https://doi.org/10.1016/j.atmosres.2025.108063

[4] Lui, G. C., Li, W. K., Leung, K. M., Lee, J. H., & Jayawardena, A. W. (2007). Modeling algal blooms using vector autoregressive model with exogenous variables and long memory filter. Ecological Modeling, 200(1-2), 130-138. https://doi.org/10.1016/j.ecolmodel.2006.06.017

[5] Qu, H. D., Liu, X., & Yang, X. P. (2026). A Legendre interpolation method for the nonlinear variable-order fractional advection–diffusion equations. Chaos, Solitons & Fractals, 204, 117770. https://doi.org/10.1016/j.chaos.2025.117770

[6] Liu, C., Chen, C., Yang, J., & Xu, Z. (2025). A fractional gray reservoir computing prediction model and its application in clean energy forecasting. Energy, 139662. https://doi.org/10.1016/j.energy.2025.139662

[7] Balagula, Y., Baimel, D., & Aharon, I. (2025). Comparing time series and neural network models of long memory for electricity price forecasting. Results in Engineering, 108465. https://doi.org/10.1016/j.rineng.2025.108465

[8] Stanisz, T., Drożdż, S., & Kwapień, J. (2024). Complex systems approach to natural language. Physics Reports, 1053, 1-84. https://doi.org/10.1016/j.physrep.2023.12.002

[9] Muravlev, A. A. (2011). Representation of a fractional Brownian motion in terms of an infinite-dimensional Ornstein-Uhlenbeck process. Russian Mathematical Surveys, 66(2), 439. https://doi.org/10.1070/RM2011v066n02ABEH004746

[10] Barndorff-Nielsen, O. E. (2001). Superposition of Ornstein–Uhlenbeck type processes. Theory of Probability & Its Applications, 45(2), 175-194. https://doi.org/10.1137/S0040585X97978166

[11] Fasen, V., & Klüppelberg, C. (2005). Extremes of supOU processes. In Stochastic analysis and applications: The Abel symposium 2005 (pp. 339-359). Springer, Berlin, Heidelberg. https://doi.org/10.1007/978-3-540-70847-6_14

[12] Siegle, P., Goychuk, I., & Hänggi, P. (2011). Markovian embedding of fractional superdiffusion. Europhysics Letters, 93(2), 20002. https://doi.org/10.1209/0295-5075/93/20002

[13] Abi Jaber, E., & El Euch, O. (2019). Multifactor approximation of rough volatility models. SIAM Journal on Financial Mathematics, 10(2), 309-349. https://doi.org/10.1137/18M1170236

[14] Di Nunno, G., & Yurchenko-Tytarenko, A. (2026). Sandwiched Volterra Volatility model: Markovian approximations and hedging. Finance and Stochastics, 30, 277-325. https://doi.org/10.1007/s00780-025-00584-2

[15] Bock, W., Desmettre, S., & da Silva, J. L. (2020). Integral representation of generalized gray Brownian motion. Stochastics, 92(4), 552-565. https://doi.org/10.1080/17442508.2019.1641093

[16] Nobis, G., Belova, A., Springenberg, M., Daems, R., Knochenhauer, C., Opper, M., ... & Samek, W. (2025). Fractional Brownian Bridges for Aligned Data. In Learning Meaningful Representations of Life (LMRL) Workshop at ICLR 2025. https://openreview.net/pdf?id=PoEFXbYNii

[17] Aichinger, F., & Desmettre, S. (2025). Pricing of geometric Asian options in the Volterra-Heston model. Review of Derivatives Research, 28(1), 5. https://doi.org/10.1007/s11147-025-09211-w

[18] Kanazawa, K., & Dechant, A. (2025). Stochastic thermodynamics for classical non-Markov jump processes. arXiv preprint arXiv:2506.04726. https://arxiv.org/abs/2506.04726

[19] Yoshioka, H., & Yoshioka, Y. (2025a). Non-Markovian superposition process model for stochastically describing concentration–discharge relationship. Chaos, Solitons & Fractals, 199, 116715. https://doi.org/10.1016/j.chaos.2025.116715

[20] Yoshioka, H., & Yoshioka, Y. (2025b). Stochastic volatility model with long memory for water quantity-quality dynamics. Chaos, Solitons & Fractals, 195, 116167. https://doi.org/10.1016/j.chaos.2025.116167

[21] Alfonsi, A. (2025). Nonnegativity preserving convolution kernels. Application to Stochastic Volterra Equations in closed convex domains and their approximation. Stochastic Processes and their Applications, 181, 104535. https://doi.org/10.1016/j.spa.2024.104535

[22] Dall'Acqua, E., Longoni, R., & Pallavicini, A. (2023). Rough-Heston local-volatility model. International Journal of Theoretical and Applied Finance, 26(06n07), 2350021. https://doi.org/10.1142/S0219024923500218

[23] Hainaut, D., Chen, J., & Scalas, E. (2025). The rough Hawkes process. Communications in Statistics-Theory and Methods, 54(11), 3322-3349. https://doi.org/10.1080/03610926.2024.2389959

[24] Lin, Y., Xu, Z., Zhang, Y., & Zhou, Q. (2025). Weighted balanced truncation method for approximating kernel functions by exponentials. Physical Review E, 112(1), 015302. https://doi.org/10.1103/xsgv-zbvp

[25] Yoshioka, H., Tanaka, T., Yoshioka, Y., & Hashiguchi, A. (2024). Statistical computation of a superposition of infinitely many Ornstein–Uhlenbeck processes. AIP Conference Proceedings, Vol. 3094, 100001. https://doi.org/10.1063/5.0210146

[26] Davis, S., Bora, B., Pavez, C., & Soto, L. (2025). Kappa distributions in the framework of superstatistics. Physica A: Statistical Mechanics and its Applications, 131191. https://doi.org/10.1016/j.physa.2025.131191
[27] Sánchez, E. (2025). Gamma-superstatistics and complex time series analysis. Physical Review E, 112(1), 014118. https://doi.org/10.1103/3kcd-3mlq
[28] Moges, E., Ruddell, B. L., Zhang, L., Driscoll, J. M., & Larsen, L. G. (2022). Strength and memory of precipitation's control over streamflow across the conterminous United States. Water Resources Research, 58(3), e2021WR030186. https://doi.org/10.1029/2021WR030186
[29] Huang, C., Zhao, X., Zhang, F., Chen, H., Song, R., Ma, G., & Cheng, W. (2025). Model averaging with logistic autoregressive conditional peak over threshold models for regional smog. Risk Analysis, 45(12), 4505-4520. https://doi.org/10.1111/risa.70069
[30] Bergmann, D. R., & Oliveira, M. A. (2026). Extreme risk clustering in long-memory financial series. Chaos, Solitons & Fractals, 202, 117513. https://doi.org/10.1016/j.chaos.2025.117513
[31] Lehrer, S., Xie, T., & Zhang, X. (2021). Social media sentiment, model uncertainty, and volatility forecasting. Economic Modelling, 102, 105556. https://doi.org/10.1016/j.econmod.2021.105556
[32] Yoshioka, H., & Yoshioka, Y. (2024a). Generalized divergences for statistical evaluation of uncertainty in long-memory processes. Chaos, Solitons & Fractals, 182, 114627. https://doi.org/10.1016/j.chaos.2024.114627
[33] Yoshioka, H., & Yoshioka, Y. (2024b). Risk assessment of river water quality using long-memory processes subject to divergence or Wasserstein uncertainty. Stochastic Environmental Research and Risk Assessment, 38(8), 3007-3030. https://doi.org/10.1007/s00477-024-02726-y
[34] Jaimungal, S., Pesenti, S. M., & Sánchez-Betancourt, L. (2024). Minimal Kullback–Leibler divergence for constrained lévy–itô processes. SIAM Journal on Control and Optimization, 62(2), 982-1005. https://doi.org/10.1137/23M1555697
[35] Kroell, E., Pesenti, S. M., & Jaimungal, S. (2024). Stressing dynamic loss models. Insurance: Mathematics and Economics, 114, 56-78. https://doi.org/10.1016/j.insmatheco.2023.11.002
[36] Kroell, E., Jaimungal, S., & Pesenti, S. M. (2025). Optimal robust reinsurance with multiple insurers. Scandinavian Actuarial Journal, 2025(5), 479-509. https://doi.org/10.1080/03461238.2024.2431539
[37] Horiguchi, S. A., & Kobayashi, T. J. (2025). Optimal control of stochastic reaction networks with entropic control cost and emergence of mode-switching strategies. PRX Life, 3(3), 033027. https://doi.org/10.1103/zttn-tpzq
[38] Yoshioka, H., & Hamagami, K. (2024). Marcus's formulation of stochastic algae population dynamics subject to power-type abrasion. International Journal of Dynamics and Control, 12(11), 3987-3999. https://doi.org/10.1007/s40435-024-01461-0
[39] Haddadchi, A., Kuczynski, A., Hoyle, J. T., Kilroy, C., Booker, D. J., & Hicks, M. (2020). Periphyton removal flows determined by sediment entrainment thresholds. Ecological Modelling, 434, 109263. https://doi.org/10.1016/j.ecolmodel.2020.109263
[40] Luce, J. J., Steele, R., & Lapointe, M. F. (2010). A physically based statistical model of sand abrasion effects on periphyton biomass. Ecological Modelling, 221(2), 353-361. https://doi.org/10.1016/j.ecolmodel.2009.09.018
[41] Pham, H. (2009). Continuous-time stochastic control and optimization with financial applications. Springer, Berlin, Heidelberg.
[42] Babonneau, F., Haurie, A., & Vielle, M. (2025). A robust infinite-horizon optimal control approach to climate economics. Central European Journal of Operations Research, 33(2), 499-528. https://doi.org/10.1007/s10100-025-00973-0
[43] Kafash, B., & Nikooeinejad, Z. (2025). A computational method for finding feedback Nash equilibrium solutions (FBNES) in Nonzero-Sum Differential Games (NZSDG) Based on the Variational Iteration Method (VIM). Mathematics and Computers in Simulation, 235, 37-59. https://doi.org/10.1016/j.matcom.2025.03.016
[44] Lykina, V., Pickenhain, S., Kolo, K., & Grass, D. (2022). Sustainability and long-term strategies in the modeling of biological processes. IFAC-PapersOnLine, 55(20), 665-670. https://doi.org/10.1016/j.ifacol.2022.09.172
[45] Marsiglio, S., & Masoudi, N. (2025). Social norms and international environmental agreements: A natural solution to environmental problems?. Energy Economics, 109049. https://doi.org/10.1016/j.eneco.2025.109049
[46] Yoshioka, H., Tsujimura, M., Hamagami, K., & Yoshioka, Y. (2020). A hybrid stochastic river environmental restoration modeling with discrete and costly observations. Optimal Control Applications and Methods, 41(6), 1964-1994. https://doi.org/10.1002/oca.2723
[47] Björk, T., Khapko, M., & Murgoci, A. (2021). Time-inconsistent control theory with finance applications. Springer, Berlin.
[48] Ekeland, I., & Lazrak, A. (2010). The golden rule when preferences are time inconsistent. Mathematics and Financial Economics, 4(1), 29-55. https://doi.org/10.1007/s11579-010-0034-x
[49] Chen, S., Luo, D., & Yao, H. (2024). Optimal investor life cycle decisions with time-inconsistent preferences. Journal of Banking & Finance, 161, 107115. https://doi.org/10.1016/j.jbankfin.2024.107115
[50] Liu, L., Niu, Y., Wang, Y., & Yang, J. (2020). Optimal consumption with time-inconsistent preferences. Economic Theory, 70(3), 785-815. https://doi.org/10.1007/s00199-019-01228-1
[51] Cetemen, D., Feng, F. Z., & Urgun, C. (2023). Renegotiation and dynamic inconsistency: Contracting with nonexponential discounting. Journal of Economic Theory, 208, 105606. https://doi.org/10.1016/j.jet.2023.105606

[52] Li, Y., Li, Z., & Zeng, Y. (2016). Equilibrium dividend strategy with nonexponential discounting in a dual model. Journal of Optimization Theory and Applications, 168(2), 699-722. https://doi.org/10.1007/s10957-015-0742-8
[53] Zhu, J., Siu, T. K., & Yang, H. (2020). Singular dividend optimization for a linear diffusion model with time-inconsistent preferences. European Journal of Operational Research, 285(1), 66-80. https://doi.org/10.1016/j.ejor.2019.04.027
[54] Kang, J. H., Huang, N. J., Yang, B. Z., & Hu, Z. (2025). Robust Equilibrium Strategy for Mean–Variance–Skewness Portfolio Selection Problem with Long Memory. Journal of Optimization Theory and Applications, 206(2), 27. https://doi.org/10.1007/s10957-025-02697-2
[55] Mbodji, O., & Pirvu, T. A. (2025). Portfolio time consistency and utility weighted discount rates. Mathematics and Financial Economics, 19, 261-291. https://doi.org/10.1007/s11579-025-00382-6
[56] Ebert, S., Wei, W., & Zhou, X. Y. (2020). Weighted discounting—on group diversity, time-inconsistency, and consequences for investment. Journal of Economic Theory, 189, 105089. https://doi.org/10.1016/j.jet.2020.105089
[57] Yoshioka, H., & Hamagami, K. (2025). Micro-macro population dynamics models of benthic algae with long-memory decay and generic growth. Results in Engineering, 32, 111783. https://doi.org/10.1016/j.rineng.2026.111783
[58] Liang, Z., & Zhang, K. (2024). Time-Inconsistent Mean Field and Agent Games under Relative Performance Criteria. SIAM Journal on Financial Mathematics, 15(4), 1047-1082. https://doi.org/10.1137/22M1533219
[59] Niu, Y., & Zou, Z. (2024). Robust dynamic contracts with multiple agents. Games and Economic Behavior, 148, 196-217. https://doi.org/10.1016/j.geb.2024.09.012
[60] Yoshioka, H., & Yoshioka, Y. (2023). Dual stochastic descriptions of streamflow dynamics under model ambiguity through a Markovian embedding. Journal of Mathematics in Industry, 13(1), 7. https://doi.org/10.1186/s13362-023-00135-4
[61] Boucekkine, R., Camacho, C., & Fabbri, G. (2013). Spatial dynamics and convergence: The spatial AK model. Journal of Economic Theory, 148(6), 2719-2736. https://doi.org/10.1016/j.jet.2013.09.013
[62] Kanazawa, K., & Sornette, D. (2024). Standard form of master equations for general non-Markovian jump processes: The Laplace-space embedding framework and asymptotic solution. Physical Review Research, 6(2), 023270. https://doi.org/10.1103/PhysRevResearch.6.023270
[63] Li, X., & Yong, J. (1995). Optimal Control Theory for Infinite Dimensional Systems. Systems & Control: Foundations & Applications. Birkhäuser, Boston.
[64] Flandoli, F., Leocata, M., Livieri, G., Morlacchi, S., Corvino, F., & Pirni, A. (2025). Structural properties in the diffusion of the solar photovoltaic in Italy: individual people/householder vs firms. Decisions in Economics and Finance. Online published. https://doi.org/10.1007/s10203-025-00532-x
[65] Yu, X., Zhang, Y., & Zhou, Z. (2021). Teamwise mean field competitions. Applied Mathematics & Optimization, 84(Suppl 1), 903-942. https://doi.org/10.1007/s00245-021-09789-1
[66] Karbach, S. (2026). Finite-rank approximation of affine processes on positive Hilbert-Schmidt operators. Journal of Mathematical Analysis and Applications, 553(2), 129852. https://doi.org/10.1016/j.jmaa.2025.129852
[67] Igwaran, A., Kayode, A. J., Moloantoa, K. M., Khetsha, Z. P., & Unuofin, J. O. (2024). Cyanobacteria harmful algae blooms: causes, impacts, and risk management. Water, Air, & Soil Pollution, 235(1), 71. https://doi.org/10.1007/s11270-023-06782-y
[68] Sun, L., Wu, L., Liu, X., Huang, W., Zhu, D., Wang, Z., ... & Liu, X. (2023). Reducing the risk of benthic algae outbreaks by regulating the flow velocity in a simulated South–North water diversion open channel. International journal of environmental research and public health, 20(4), 3564. https://doi.org/10.3390/ijerph20043564
[69] Xu, C., Jia, T., Xu, T., Lan, Y., Li, N., & Jia, H. (2025). Deep learning-based optimal adaptive regulation pathway of algal blooms in urban rivers under long-term uncertainties. Water Research, 124677. https://doi.org/10.1016/j.watres.2025.124677
[70] Yang, M., Guan, G., Bi, Y., Zhu, Y., & Wang, K. (2025). A coupled hydraulic-ecological model for simulating periphytic algal detachment in water delivery canals. Journal of Environmental Management, 384, 125543. https://doi.org/10.1016/j.jenvman.2025.125543
[71] Graba, M., Moulin, F. Y., Boulêtreau, S., Garabétian, F., Kettab, A., Eiff, O., ... & Sauvage, S. (2010). Effect of near-bed turbulence on chronic detachment of epilithic biofilm: Experimental and modeling approaches. Water Resources Research, 46(11). https://doi.org/10.1029/2009WR008679
[72] Hamagami, K., Yoshioka, H., & Ito, J. (2024). Effect of Bed Forms on the Detachment of Benthic Algae by the Sediment Transportation. Journal of Water and Environment Technology, 22(4), 159-167. https://doi.org/10.2965/jwet.23-127
[73] Clason, C. (2020). Introduction to Functional Analysis. Birkhäuser, Cham.
[74] Yoshioka, H., & Yoshioka, Y. (2024b). Statistical evaluation of a long-memory process using the generalized entropic value-at-risk. Environmetrics, 35(4), e2838. https://doi.org/10.1002/env.2838
[75] Maurya, J., & Misra, A. K. (2025). Modeling the impact of harmful algal blooms on aquatic life and human health. The European Physical Journal Plus, 140(5), 391. https://doi.org/10.1140/epjp/s13360-025-06295-z
[76] Sajid, M., Misra, A. K., & Almohaimeed, A. S. (2024). Modeling the role of fish population in mitigating algal bloom. Electronic Research Archive, 32(10). https://doi.org/10.3934/era.2024269

[77] Beal, M. R., & Schaeffer, B. (2026). Pixel-scale satellite forecasting of cyanobacteria in Florida lakes. Harmful Algae, 103041. https://doi.org/10.1016/j.hal.2025.103041
[78] Yoshioka, H. (2024). Modeling stationary, periodic, and long memory processes by superposed jump-driven processes, Chaos, Solitons & Fractals, 188, 115357. https://doi.org/10.1016/j.chaos.2024.115357